\documentclass[12pt]{amsart}
\usepackage{hyperref}
\usepackage{latexsym}
\usepackage{amsmath}
\usepackage{amsfonts}
\usepackage{amssymb}
\usepackage{color}
\usepackage{ulem}
  \newcommand{\B}{\mathcal{B}}
 \newcommand{\cB}{\mathcal{B}}

\newcommand{\cD}{\mathcal{D}}

\newcommand{\G}{\mathcal{G}}
\newcommand{\cZ}{\mathcal{Z}}
\newcommand{\gG}{\Gamma}

\newcommand{\T}{\mathbb{T}}
\newcommand{\Q}{\mathbb{Q}}
\newcommand{\R}{\mathbb{R}}
\newcommand{\E}{\mathbb{E}}
\newcommand{\N}{\mathbb{N}}
\newcommand{\bbZ}{\mathbb{Z}}
\newcommand{\Z}{\mathbb{Z}}

\newcommand{\norm}[1]{\left\Vert #1\right\Vert}
\newcommand{\nnorm}[1]{\lvert\!|\!| #1|\!|\!\rvert}

\theoremstyle{plain}
\newtheorem{theorem}{Theorem}[section]
\newtheorem{lemma}[theorem]{Lemma}
\newtheorem{proposition}[theorem]{Proposition}

\newtheorem*{conjecture1}{Conjecture I}
\newtheorem*{conjecture2}{Conjecture II}

\newtheorem{bigthm}{Theorem}

\newtheorem{bigthmp}{Theorem}

\newtheorem{question}{Question}
\newtheorem*{theorem*}{Theorem}
\newtheorem*{Correspondence1}{Furstenberg Correspondence Principle}

\theoremstyle{definition}
\newtheorem{definition}[theorem]{Definition}

\newtheorem*{example*}{Example}

\theoremstyle{remark}

\newtheorem*{remark}{Remark}
\newtheorem*{remarks}{Remarks}

\begin{document}
\title{Powers of sequences and recurrence}
\author{Nikos Frantzikinakis}
\address[Nikos  Frantzikinakis]{Department of Mathematics\\
 University of Memphis\\
  Memphis, TN \\ 38152 \\ USA } \email{frantzikinakis@gmail.com}

\author{Emmanuel Lesigne}
\address[Emmanuel Lesigne]{Laboratoire de Math\'ematiques et Physique Th\'eorique (UMR CNRS 6083)\\
Universit\'e Fran\c{c}ois Rabelais Tours \\
F\'ed\'eration de Recherche Denis Poisson\\
Parc de Grandmont \\
37200 Tours\\ France} \email{emmanuel.lesigne@lmpt.univ-tours.fr}

\author{M\'at\'e Wierdl}
\address[M\'at\'e Wierdl]{Department of Mathematics\\
  University of Memphis\\
  Memphis, TN \\ 38152 \\ USA } \email{wierdlmate@gmail.com}

\begin{abstract}
We study recurrence, and multiple recurrence, properties along the
$k$-th powers of a given set of integers. We show that the
property of recurrence for some given values of $k$ does not give
any constraint on the recurrence for the other powers. This is
motivated by similar results in  number theory concerning additive
basis of natural numbers. Moreover, motivated by a result of Kamae
and Mend\`es-France, that links single recurrence with uniform
distribution properties of sequences, we look for an analogous
result dealing with higher order  recurrence  and make a related
conjecture.
\end{abstract}

\thanks{The  first author was partially supported
 by NSF grant DMS-0701027.}



\maketitle

\setcounter{tocdepth}{1}
\tableofcontents
\section{Motivation, historical remarks}
\label{sec:motiv-hist-remarks}

In their 1912 paper (\cite{HL}), Hardy and Littlewood
proved that for any irrational number  $\alpha$, the set
$(\alpha,\alpha),(2\alpha,2^2\alpha),(3\alpha,3^2\alpha),\dots$ $\mod
1$ is dense in the two dimensional torus, and the obvious extension of
this result to higher dimensional tori, and arbitrary powers.  This result can be considered as an indication that in
some sense the sequences $1,2,3\dots$ and $1^2,2^2,3^2\dots$ behave
{\it independently}.  Soon after Hardy-Littlewood's work, Weyl proved
that, in fact, the sequence
$(\alpha,\alpha),(2\alpha,2^2\alpha),(3\alpha,3^2\alpha),\dots$ $\mod
1$ is uniformly distributed in the two dimensional torus (\cite{We}).

Since these results, a number of theorems appeared that could be
interpreted to express the independence of sequences of powers of $n$.
A direct motivation for our work is the result of Deshouillers, Erd\"os
and S\'ark\"ozy  (\cite{DES}) on bases.  A set $B$ of positive
integers is called a {\it basis} if there exists an $h\in\N$  such that every
 $n\in \N$ can be written in the form $n=b_1+b_2+\dots+b_h$, $b_i\in
B\cup \{0\}$.  The result of Deshouillers, Erd\"os and S\'ark\"ozy, says that
there exists a set of integers $B=\{b_1,b_2,\ldots \}$ which is not a basis, but
the set $B^2=\{b_1^2,b_2^2,\ldots \}$  {\it is} a
basis.  They also construct a set of integers which is a basis, but
the set  of squares of its elements is not a basis.  The ultimate generalization of this
result appears in a paper by Deshouillers and Fouvry
(\cite{DF}).  It says the following: given any set $G$  (the set of  ``good" exponents), there exists a set of integers
$B=\{b_1,b_2,\ldots \}$ such  that, for  $k\in\N$, the set $B^k=\{b_1^k,b_2^k,\ldots \}$ is a basis if and only if
$k\in G$.  Again, this result can be interpreted as expressing the
independent behavior of powers, and indeed, the proof utilizes a
quantitative version of Weyl's theorem via the Hardy-Littlewood circle
method.

In this paper we deal with  intersective sets. Recall that a set of integers has {\it positive upper density} if $\bar{d}(\Lambda)=\limsup_{N\to\infty}|\Lambda\cap
\{-N,\ldots ,N\}|/(2N+1)>0$ (if the $\limsup$ is a limit we denote it by $d(\Lambda)$). A set of integers $R$ is
called {\it intersective},  if for every set of integers  $\Lambda$  with positive
upper density, the equation $x-y=r$ is solvable in $x,y\in \Lambda$ and $r\in
R\setminus{ \{0\}}$. The adjective ``intersective'' was  introduced by Ruzsa (\cite{Ru}),
because the
solvability of the equation $x-y=r$ with $x,y\in \Lambda$ and $r\in R\setminus{ \{0\}}$ can
be expressed in the following way:
\begin{equation*}
  \Lambda \cap (\Lambda-r)\neq \emptyset \ \text{  for some nonzero } r\in R.
\end{equation*}
It is a well known theorem of S\'ark\"ozy (\cite{Sa}), that for
any fixed $k$, the set of $k$-th powers is intersective.
Bergelson and H\aa land (\cite{BH})
constructed a set of integers that is not intersective but the
set of squares of its elements is intersective. We establish
a generalization of their result to the setting of the
Deshouillers--Fouvry theorem we mentioned above:
\begin{bigthm}
  \label{T:A}
  Let $G$ be a set of positive integers.

There exists a set of integers $R=\{r_1,r_2,\ldots \}$ such that: for  $k\in \N$, the set
  $R^k=\{r_1^k,r_2^k,\ldots \}$  is intersective if and only $k\in G$.
\end{bigthm}
In our paper we go further, and prove a generalization of this result
to arithmetic progressions.  The famous result of Szemer\'edi (\cite{Sz}) states that
every set  of integers with positive upper density contains arbitrarily long
arithmetic progressions.  Bergelson and Leibman proved in
\cite{BL} that, for any fixed $k$, every set of positive upper
density contains arbitrarily long arithmetic progressions so that the
difference of the progression is a $k$-th power.  To be able to talk
about our result, we introduce the following definition.  A set of integers $R$ is
called  {\it $\ell$-intersective} if every set $\Lambda$ of positive
upper density contains an $\ell+1$ long
arithmetic progression with a common difference from $R\setminus \{0\}$.  In other words, $\Lambda$
contains configurations of the form
$$m, m+r, m+2r,\ldots, m+\ell r$$
for some nonzero $r\in R$.  Similarly to single
intersectivity, we can rewrite this condition  as
\begin{equation*}
  \Lambda \cap (\Lambda-r) \cap(\Lambda-2r)\cap\dots \cap
  (\Lambda-\ell r)\ne \emptyset \  \text{ for some  nonzero } r\in R.
\end{equation*}
In this language our
result is the exact ``multiple'' analog of
Theorem~\ref{T:A}:
\begin{bigthm}
  \label{T:B}
  Let $\ell$ be a positive integer and  $G$ be a set of
  positive integers.

  There exists a set of integers $R=\{r_1,r_2,\ldots \}$  such that: for  $k\in \N$, the set
  $R^k=\{r_1^k,r_2^k,\ldots\}$ is $\ell$-intersective if and only $k\in G$.
\end{bigthm}
Note that, for example,  we do not claim the existence of a set $R$  that is  $\ell$-intersective {\it for every} $\ell\in \N$, but the set $R^2$ is not $\ell$-intersective for some  $\ell\in \N$. Actually, it may very well not be
 possible to construct such an example (see Question~\ref{Q:question2} in Section~\ref{SS:questions}).

We will use Furstenberg's correspondence principle  to translate the previous statements to ergodic theory
and  will then  verify the corresponding ergodic statements.
In the next section we explain the ergodic theoretical analogs of the previous
theorems, and we state our main results.

\section{Main results}
\label{sec:main-results-conj}
\subsection{Good and bad powers for sets of recurrence and  $\ell$-recurrence}
All along the article we will use the word  {\it system}, or the term {\it measure preserving system},
 to designate a quadruple
$(X,\mathcal{B},\mu, T)$, where $(X,\mathcal{B},\mu)$ is a
probability space, and $T\colon X\to X$ is an \emph{invertible}   measurable map such
that $\mu(T^{-1}A)=\mu(A)$ for all $A\in\mathcal{B}$.
 In \cite{Fu1}, Furstenberg perceived a
connection between existence of structures in sets of integers having
positive upper density and recurrence properties of measure preserving
systems. He used this, and other  ideas, to  give an ergodic theoretic proof of
Szemer\'edi's theorem on arithmetic progressions. This new approach gave rise to the field of ergodic
Ramsey theory, where problems in combinatorial number theory
are treated using techniques from ergodic theory, and led to several far-reaching extensions of Szemer\'edi's theorem.
 We will use this correspondence to translate statements about  ``intersectivity''
 to statements about  ``recurrence''.
 The
following formulation is from \cite{B}:

\begin{Correspondence1}\label{correspondence} Let $\Lambda$ be a
set of  integers. There exist a system $(X,\B,\mu,T)$ and a set $A\in\mathcal{B}$ with
$\mu(A)=\bar{d}(\Lambda)$ such that
\begin{equation}\label{E:correspondence}
\bar{d}(\Lambda \cap (\Lambda-n_1)\cap\ldots\cap
(\Lambda-n_\ell))\geq \mu(A\cap T^{-n_1}A\cap\cdots \cap T^{-n_\ell}A),
\end{equation}
for all  $n_1,\ldots,n_\ell\in\Z$ and $l\in\N$.
\end{Correspondence1}
Using this principle, we can  reformulate
Theorems~\ref{T:A} and \ref{T:B} in ergodic theoretic language. We first translate the notion of
 intersectivity.
\begin{definition}
  We say that the set of integers $R$ is a  {\it set of recurrence for
    the  system $(X,\mathcal{B},\mu,T)$}, if for every set $A\in
  \mathcal{B}$ of positive measure, we have
  \begin{equation}\label{E:one}
    \mu(A\cap T^{-r}A)>0 \text{ for some nonzero } r\in R.
  \end{equation}
  We say that the set of integers $R$  is   {\it a set of recurrence}, or {\it good for recurrence}, if
  it is a set of recurrence for every system.
\end{definition}
Note that  if $R$ is a set of recurrence then \eqref{E:one} is in fact
satisfied for infinitely many $r\in R$.

Using Furstenberg's correspondence principle it is easy to see that if  $R$ is a set of
recurrence then it is intersective (the converse is   also true and not hard to show).
As a consequence, the following result implies
Theorem~\ref{T:A}:
\begin{bigthmp}
  \label{T:A'}
  Let $G$ be a set of positive integers.

  There exists a set of integers $R=\{r_1,r_2,\ldots \}$   such that: for  $k\in \N$,
   the set $R^k=\{r_1^k,r_2^k,\ldots \}$  is good for recurrence if and only $k\in G$.
\end{bigthmp}
Although Theorem~\ref{T:A'} will be later subsumed by a stronger result (Theorem~\ref{T:B'}), we choose
to give an independent proof of it in Section~\ref{S:single}, as in this case the analysis does not depend upon
complicated multiple ergodic theorems, and so it becomes easier to
see the main ideas of the proof.

Now to formulate the ergodic theoretical analog of
Theorem~\ref{T:B}, we first translate the notion of  $\ell$-intersectivity.
\begin{definition}
  Let $\ell$ be a positive integer.  We say that the set $R$ of
  integers is   {\it a set of $\ell$-recurrence for the system
    $(X,\mathcal{B},\mu,T)$}, if for every set $A\in \mathcal{B}$ of
  positive measure, we have
  \begin{equation*}
    \mu(A\cap T^{-r}A \cap T^{-2r}A\cap\dots \cap T^{-\ell r}A )>0
    \text{ for some nonzero } r\in R.
  \end{equation*}
  We say that the set of integers $R$  is  {\it a set of $\ell$-recurrence}, or {\it  good for  $\ell$-recurrence}, if
  it is a set of $\ell$-recurrence for every system.
\end{definition}
Using Furstenberg's correspondence principle it is easy to see that if  $R$ is a set of
$\ell$-recurrence then it is $\ell$-intersective. The converse is also true and easy to establish (\cite{BHRF}).
Some examples of sets of
$\ell$-recurrence (or $\ell$-intersective), for every $\ell\in \N$, are IP sets, meaning sets that consist of all finite sums (with distinct entries)
of some infinite set (\cite{FuK2}),
and sets of the form $\bigcup_{n\in\N}\{a_n,2a_n,\ldots,na_n\}$
where $(a_n)$ is a sequence of nonzero integers (this follows from a finite version of Szemer\'edi's theorem). It is also known that the
set of values of any non-constant integer polynomial with
zero constant term is a set of $\ell$-recurrence for all $\ell\in\N$ (\cite{BL}).
Examples of sets which are not sets of recurrence are sets which do not contain any multiple of a given $d\in\N$, and lacunary sets. An example of a set of $\ell$-recurrence but not $(\ell+1)$-recurrence
is $\big\{n\in\N\colon \{n^{\ell+1}\alpha\}\in [1/4,3/4] \big\}$ where $\alpha$ is any irrational number (\cite{FLW}).

The ergodic theoretical analog of
Theorem~\ref{T:B} is:
\begin{bigthmp}
  \label{T:B'}
  Let $\ell$ be a positive integer and  $G$ be a set of positive integers.

  There exists a set of integers  $R=\{r_1,r_2,\ldots \}$   such that: for  $k\in \N$, the set
   $R^k=\{r_1^k,r_2^k,\ldots \}$  is good for  $\ell$-recurrence if and only $k\in G$.
\end{bigthmp}
We prove this result in
Section~\ref{S:multiple1} using an argument  similar to the one used to prove Theorem~\ref{T:A'}.
 There are some extra difficulties in this case though, since we have to establish   some equidistribution results on nilmanifolds, and also a  uniform multiple recurrence result
that we state and prove in the Appendix.
\subsection{Good and bad powers for  sets of $(n^{a_1},\ldots,n^{a_\ell})$-recurrence}
The previous results   deal with sets of recurrence along families of
polynomials of the form $\{n^k,2n^k,\ldots,\ell n^k\}$. We also study
the other extreme case, of sets of recurrence along families of
linearly independent polynomials, like families of the form
$\{n^{a_1},\ldots, n^{a_{\ell}}\}$, where $a_1,\ldots,a_{\ell}\in \N$ are
distinct. The next definition will facilitate our discussion.
\begin{definition}
Let $u_1(n),\ldots,u_{\ell}(n)$ be  integer sequences. We say that the
set $R\subset \N$ is {\it good for recurrence along the sequence
$(u_1(n),\ldots,u_{\ell}(n))$}, if for every system $(X,\mathcal{B},\mu,T)$ and $A\in \mathcal{B}$
of positive measure, there exist infinitely many $r\in R$ such that
$$
\mu (A\cap T^{-u_1(r)}A\cap \ldots\cap T^{-u_{\ell}(r)}A)>0.
$$
\end{definition}
Using Furstenberg's correspondence principle it is easy to define an analogous notion in combinatorics.
We will show:
\begin{bigthm}\label{T:C}
Let $A_{\ell}=\{(a_1,a_2,\ldots,a_{\ell})\in \mathbb{N}^{\ell}\colon
a_1<a_2<\ldots<a_{\ell}\}$ and  $G\subset A_{\ell}$.

There exists a set of integers $R$ such that:
for $ (a_1,a_2,\ldots,a_\ell) \in A_{\ell}$ the set 
$R$ is good for recurrence along the sequence
$(n^{a_1},\ldots,n^{a_{\ell}})$ if and only if $(a_1,\ldots, a_{\ell})\in
G$.
\end{bigthm}
The proof of this result is  similar to the proof of Theorem~\ref{T:B'} and we give it in Section~\ref{S:multiple2}.

\subsection{Powers of sequences and sufficient conditions for  $\ell$-recurrence}
When $\ell\geq 2$, there is currently  no general criterion providing usable sufficient conditions
 for a set of positive integers $R=\{r_1,r_2,\ldots \}$  to be good for $\ell$-recurrence.
In contrast, when $\ell=1$ such a criterion exists, it is a
 result of Kamae and
Mend\`es-France~\cite{KM}  that links recurrence properties of a set  $R$ with uniform distribution
properties of   sequences of the form $(r_n\alpha)_{n\in \N}$  in $\T$, where $\alpha$ is irrational  (see Theorem~\ref{T:single}).

In Section~\ref{S:conjecture} we
are looking for  a similar result for
$\ell$-recurrence. When $\ell=2$ it is well understood that such a criterion should
be related to  stronger uniform distribution properties of ``quadratic
nature", forcing  at the very least good single recurrence properties for the sequence $(r_n)_{n\in \N}$ and the sequence of squares
$(r_n^2)_{n\in \N}$. We show why a plausible statement involving only quadratic functions of $r_n$ fails (Theorem~\ref{T:main2}), the reason
turns out to be that
one has to also  take into account generalized quadratic functions of $r_n$, that is, sequences like $([r_n\alpha]r_n)_{n\in\N}$ where $\alpha$ is irrational (Lemma~\ref{L:1}).
 Using the language of nilsystems we state a conjecture  that, if true,  would provide a natural generalization of the Kamae and Mend\`es-France criterion for $\ell$-recurrence.
 In Theorem~\ref{T:conjecture}  we verify this conjecture when the set $R$ has positive density.

\bigskip
{\bf Notation:} The following notation will be used throughout the
article: $Tf=f\circ T$, $e(x)=e^{2\pi i x}$, $\{x\}=x-[x]$. For a bounded numerical sequence $(a_n)_{n\in\N}$ we will write
$\text{D-}\!\lim_{n\to\infty}(a_n)=0$ if one of the following equivalent properties is satisfied :\\
- For every $\varepsilon>0$, $d \big( \{ n: |a_n|>\varepsilon\} \big) =0$;\\
- There exists $E\in\N$, with $d(E)=1$, such that $\lim_{n\to\infty,n\in E}\ a_n=0$;\\
- $\lim_{N\to+\infty}\frac1N\sum_{n=1}^N |a_n| =0$.

\section{Good and bad powers for sets of  recurrence}\label{S:single}
In this section we will prove Theorem~\ref{T:A'}.
The proof is based on the following ergodic  result:
\begin{proposition}\label{P:ergodic1}
Let $(X,\mathcal{B},\mu, T)$ be a system, $f\in L^\infty(\mu)$, $h_1,\ldots,h_{s}\colon
\mathbb{T}\to\mathbb{C}$ be Riemann integrable functions, and
$\beta$ be an irrational number. If $k,k_1,\ldots,k_s$ are
distinct positive integers  then
\begin{equation}\label{E:11}
\lim_{N\to\infty}\frac{1}{N}\sum_{n=1}^N
h_1(n^{k_1}\beta)\cdot\ldots\cdot h_s(n^{k_s}\beta)\cdot
T^{n^{k}}f=\int h_1 \ dt \cdot\ldots \cdot \int h_s \ dt\cdot
 \lim_{N\to\infty}\frac{1}{N}\sum_{n=1}^N T^{n^k} f,
\end{equation}
where the convergence takes place in $L^2(\mu)$.
\end{proposition}
\begin{remark}
It is well known, as a direct consequence of the spectral theorem
and Weyl's uniform distribution theorem,
 that both limits exist in $L^{2}(\mu)$.
\end{remark}
\begin{proof}
 Using a standard estimation by
continuous functions from above and
  below (for every $\varepsilon>0$ we can
  find continuous functions $\underline{h}, \overline{h}$,
  such that $\underline{h}\leq h\leq  \overline{h}$ and
  $\int(\overline{h}- \underline{h}) \ dt \leq \varepsilon$),
it suffices to check that \eqref{E:11} holds
   when $h_1,\ldots,h_s$ are  continuous
  functions. Using Weierstrass approximation theorem of continuous
  functions by trigonometric polynomials, and  linearity, it
  suffices to check the result when
  $h_1(t)= e(l_1 t),\ldots, h_s(t)= e(l_s t)$ for some
  $l_1,\ldots,l_s \in\mathbb{Z}$ (remember $e(x)=e^{2\pi i x}$). If all the $l_i$'s are zero
  then \eqref{E:11} holds trivially. So without loss of generality we can assume that
  $l_1\neq 0$.
Using the spectral theorem for the unitary action $T$ on $L^2(\mu)$, we associate to the function $f$ a finite positive measure $\sigma_f$ on the torus $\T$, such that, for any complex numbers $a_1,\ldots, a_N$ we have
\begin{equation}\label{E:spectral}
\norm{\frac{1}{N}\sum_{n=1}^N a_n \cdot T^{n^k}f}_{L^2(\mu)}=
\norm{\frac{1}{N}\sum_{n=1}^N a_n \cdot  e(n^k t)}_{L^2(\sigma_f(t))}.
\end{equation} Setting $a_n=e(l_1n^{k_1}\beta)\cdot\ldots\cdot
e(l_sn^{k_s}\beta)$ in \eqref{E:spectral}, we see that it suffices
to show that
\begin{equation}\label{E:spectral2}
\lim_{N\to\infty}\norm{\frac{1}{N}\sum_{n=1}^N
e(l_1n^{k_1}\beta)\cdot\ldots\cdot e(l_sn^{k_s}\beta) \cdot  e(n^k
t)}_{L^2(\sigma_f(t))}=0.
\end{equation}
  The average in
  \eqref{E:spectral2} can be written as
  $$
\frac{1}{N}\sum_{n=1}^N e\big(l_1n^{k_1}\beta+\ldots+l_s
n^{k_s}\beta+n^kt \big).
  $$
Since the integers $k, k_1,\ldots,k_s$  are distinct, the
coefficient of $n^{k_1}$ is $l_1\beta$, which  is irrational since
$\beta$ is irrational and $l_1\neq 0$. By Weyl's uniform
distribution theorem the last average converges
 to zero for every $t\in [0,1)$,
  which gives \eqref{E:spectral2}.
\end{proof}

We will also use the  following result which was proved in \cite{Fo} and \cite{BH} (a
more general result is proved in the Appendix of the present paper):
\begin{theorem}\label{T:uniform}
Suppose that  $R$ is a set of single recurrence. For every
$\varepsilon>0$ there exists an $N=N(\varepsilon)$ and
$\delta=\delta(\varepsilon)>0$, such that for every  system
$(X,\mathcal{B},\mu,T)$ and set $A\in \mathcal{B}$ with
$\mu(A)\geq \varepsilon$ we have
$
 \mu(A\cap T^{-n}A)\geq \delta
$
for some nonzero $n\in R\cap[-N,N]$.
\end{theorem}
\begin{proof}[Proof of Theorem~\ref{T:A'}]

We first consider the case where the complement $B$ of $G$ is
empty or finite. If $B$ is empty then $R=\mathbb{N}$ works (see \cite{Fu2}).
If $B$ is finite then  $B=\{b_1,\ldots ,b_s\}$ for some $b_i\in\N$.
Fix an irrational number $\beta$, and for $k\in \N$ let
$$
R_k=\{n\in\mathbb{N}\colon \{n^k\beta\}\in [1/4,3/4]\}.
$$ We claim
that the set $R= R_{b_1}\cap\ldots \cap R_{b_s}$ has the advertised
property.

The  set $R^b$ is not good for single recurrence for $b\in B$
since it is not good for recurrence for the rotation by $\beta$ on
$\mathbb{T}$. Now take $g\in G$. Using
Proposition~\ref{P:ergodic1} we will show that $R^{g}$ is a set of
single recurrence. Let $(X,\mathcal{B},\mu, T)$ be a system and
$A\in \mathcal{B}$ with $\mu(A)>0$. Set $h_i={\bf
1}_{[1/4,3/4]}$, $k_i=b_i$,  for $i=1,\ldots,s$, $k=g$,  and $f={\bf 1}_A$ in
\eqref{E:11}, then multiply by $ {\bf 1}_A$ and integrate with
respect to $\mu$. We get
$$
\lim_{N\to\infty} \frac{1}{N}\sum_{1\leq n\leq N, n\in R}
\mu(A\cap T^{-n^g}A)= \frac{1}{2^s}
\lim_{N\to\infty}\frac{1}{N}\sum_{n=1}^N \mu(A\cap T^{-n^g}A).
$$
 The last limit is positive (this is implicit in  \cite{Fu1} and
 \cite{Fu2} and explicit in \cite{BL}), showing that $R^g$
 is a set of single recurrence.

Now  we deal with the general case where the set of bad powers $B$ is  infinite.
 Fix $s\in \N$ and let $\tilde{R}_s=
\bigcap_{b\in B, b\leq s} R_b$. We showed before that $(\tilde{R}_s)^g$ is a set
of single recurrence for $g\in G$. By Theorem~\ref{T:uniform}
 there exists a finite set
$F_s\subset \tilde{R}_s$ such that for each  $g\in G \cap [1,s]$ the following
is true: For every  system $(X,\mathcal{B},\mu,T)$ and every
$A\in\mathcal{B}$ with $\mu(A)>1/s$ we have that $\mu(A\cap
T^{-n}A)>0$ for some $n\in (F_s)^g$.

We claim that  $R=\bigcup_{s\in\N} F_s$ is the set we are looking
for. We first show that $R^b$ is not a set of single recurrence
for $b\in B$. Let  $b$ be an element of $B$. Since $R$ is
contained in $R_b$ up to a finite set,  and $(R_b)^{b}$ is not
a set of single recurrence,  we conclude that $R^{b}$ is not a set
of single recurrence. Suppose now that $g\in G$, it remains to
show that $R^g$ is a set of single recurrence. Let
$(X,\mathcal{B},\mu,T)$ be a system and $A\in\mathcal{B}$ be such
that $\mu(A)>0$. Then $\mu(A)>1/s$ for some $s\in \N$ with $s>g$.
By the definition of $F_s$ we have that $
 \mu(A\cap T^{-n}A)>0
$ for some $n\in (F_s)^g$. Since $F_s\subset R$ we conclude that
$R^g$ is a set of single recurrence. This completes the proof.
\end{proof}

\section{Background in ergodic theory } \label{SS:ergodic}
\subsection{Factors in ergodic theory}
 Throughout the article we consider {\it invertible} measure preserving systems $(X,\mathcal{B},\mu, T)$ where the probability space $(X,\mathcal{B},\mu)$ is a {\it Lebesgue space}. This classical assumption allows us to use Rokhlin's theory of factors and disintegration.
(The basic reference here is \cite{Ro}, see also \cite{Zi} Section 1.1, \cite{Rud} Chapter 2, or \cite{W} Section~2.3.) These two extra assumptions are not at all restrictive for our purposes, the reason being that  the measure preserving systems constructed using the correspondence principle of Furstenberg are invertible and Lebesgue.

A {\it homomorphism} from a system $(X,\mathcal{B},\mu, T)$ onto
a system $(Y, \cD, \nu, S)$
 is a measurable map $\pi\colon X'\to Y'$, where $X'$ is a $T$-invariant subset
 of $X$ and $Y'$
 is an $S$-invariant subset of $Y$, both of full measure, such that $\mu\circ\pi^{-1} = \nu$
 and $S\circ\pi(x) = \pi\circ T(x)$
  for $x\in X'$. When we have such a homomorphism we say
  that the system $(Y, \cD, \nu, S)$ is a {\it factor} of the
  system $(X,\mathcal{B},\mu, T)$.    If the  factor map $\pi\colon X'\to Y'$
can be chosen to be  injective, then we say that the systems
$(X,\cB, \mu, T)$ and $(Y, \cD, \nu, S)$ are {\it isomorphic}
(bijective maps on Lebesgue spaces have measurable inverses).


A factor can be characterized (modulo isomorphism) by the data
$\pi^{-1}(\mathcal{D})$ which is a $T$-invariant
sub-$\sigma$-algebra of $\mathcal B$, and any $T$-invariant
sub-$\sigma$-algebra of $\mathcal B$ defines a factor; by a
classical abuse of terminology we denote by the same letter the
$\sigma$-algebra $\mathcal{D}$ and its inverse image by $\pi$. In other
words, if $(Y, \cD, \nu, S)$ is a factor of
$(X,\mathcal{B},\mu,T)$, we think of $\cD$ as a sub-$\sigma$-algebra of
$\mathcal{B}$. A factor can also be characterized (modulo
isomorphism) by a $T$-invariant sub-algebra $\mathcal{F}$ of
 $L^\infty(X,\mathcal{B},\mu)$, in which case $\cD$ is the sub-$\sigma$-algebra
 generated by $\mathcal{F}$, or
 equivalently, $L^2(X,\cD,\mu)$ is the closure of $\mathcal{F}$ in $L^2(X,\mathcal{B},\mu)$.
 We will sometimes abuse notation and use the sub-$\sigma$-algebra $\mathcal{D}$ in place of
 the sub-space $L^2(X,\cD,\mu)$. For example, if we write that a function is orthogonal to the factor $\mathcal{D}$, we mean that  is orthogonal to the sub-space $L^2(X,\cD,\mu)$.

If $\cD$ is a $T$-invariant sub-$\sigma$-algebra of $\cB$ and
$f\in L^2(\mu)$, we define the {\it conditional expectation
$\mathbb{E}(f|\cD)$ of $f$ with respect to $\cD$} to be the
orthogonal projection of $f$ onto $L^2(\cD)$. We frequently make use
of the identities
$$
\int \mathbb{E}(f|\cD) \ d\mu= \int f\ d\mu, \quad
T\,\mathbb{E}(f|\cD)=\mathbb{E}(Tf|\cD).
$$
(If we want to indicate  the dependence on the reference measure, we write $\mathbb{E}=\mathbb{E}_\mu$.)

For each $d\in\mathbb{N}$, we define $\mathcal{K}_d$ to be the
factor induced by the function algebra
$$\{f\in L^\infty(\mu):T^df=f\}.$$
We define the {\it rational Kronecker factor} $\mathcal{K}_{rat}$
to be the factor induced by the algebra generated by the functions
$$
\{f\in L^\infty(\mu):T^df=f \text{ for some } d\in\mathbb{N}\}\ .
$$
This algebra is the same as the
algebra spanned by the
 bounded functions that satisfy $Tf=e(a)\cdot
f$ for some $a\in \Q$.

The {\it Kronecker factor} $\mathcal{K}$  is induced by the
algebra spanned by the
 bounded
eigenfunctions of $T$, that means, functions that satisfy $Tf=e(a)\cdot
f$ for some $a\in \R$.

It is known that if $f$ is a bounded function such that $\mathbb{E}_\mu(f|\mathcal{K}(T))=0$, then $\mathbb{E}_{\mu\otimes\mu}(f\otimes f|\mathcal{K}_{rat}(T\times T))=0$ (see for example \cite{Fu2}, Section 4.4).

The transformation $T$ is {\it ergodic} if $Tf=f$ implies that
$f=c$ (a.e.) for some $c\in \mathbb{C}$.
Every system $(X,\cB,\mu,T)$ has an {\it ergodic decomposition},
meaning that we can write $\mu=\int \mu_t\ d\sigma(t)$, where
$\sigma$ is a probability measure on $[0,1]$ and $\mu_t$ are
$T$-invariant probability measures on $(X,\cB)$ such that the
systems $(X,\cB,\mu_t,T)$ are ergodic for $t\in [0,1]$. We
sometimes denote the ergodic components by $T_t, t\in [0,1]$.

We say that $(X,\cB,\mu,T)$ is an {\it inverse limit of a
sequence of  factors} $(X,\cB_j,\mu,T)$ if
$(\cB_j)_{j\in\mathbb{N}}$ is an increasing sequence of
$T$-invariant sub-$\sigma$-algebras such that
$\bigvee_{j\in\N}\mathcal{B}_j=\mathcal{B}$ up to sets of measure
zero.

\subsection{Characteristic factors}
Following \cite{HK1}, for every  system $(X,\mathcal{B},\mu,T)$
and function $f\in L^\infty(\mu)$, we define inductively the
(functional valued) seminorms $\nnorm{f}_\ell$ as follows: For $\ell=1$
we set $\nnorm{f}_1=|\E(f|\mathcal{I})|$, where $\mathcal{I}$ is
the $\sigma$-algebra of $T$-invariant sets. For $\ell\geq 2$ we set
\begin{equation}
\label{eq:recur} \nnorm f_{\ell+1}^{2^{\ell+1}} =\lim_{N\to+\infty}\frac
1N\sum_{n=1}^{N} \nnorm{\overline{f}\cdot T^nf}_{\ell}^{2^{\ell}}.
\end{equation}
It was shown in~\cite{HK1}  that for
every integer $\ell\geq 1$, $\nnorm\cdot_\ell$ is a seminorm on
$L^\infty(\mu)$ and  it defines factors $\cZ_{\ell-1}=\cZ_{\ell-1}(T)$
in the following manner: the $T$-invariant sub-$\sigma$-algebra
$\cZ_{\ell-1}$ is characterized by
$$
\text{ for } f\in L^\infty(\mu),\  \E(f|\cZ_{\ell-1})=0\text{ if and
only if } \nnorm f_{\ell} = 0.
$$
(In \cite{HK1} the authors
work with ergodic systems, in which case $\nnorm{f}_1=\left|\int f
\ d\mu\right|$,  and real valued functions, but the whole discussion can be carried out for
nonergodic systems  and complex valued functions without extra difficulties.)
 If $f$ is a bounded function such that $\E_\mu(f|\cZ_\ell(T))=0$ then
 $\E_{\mu\otimes\mu}(f\otimes f|\cZ_{\ell-1}(T\times T))=0$ (this is implicit in
\cite{HK1}). Also, if $T_t$ where  $t\in [0,1]$ are the ergodic
components of the system, then $\E(f|\mathcal{Z}_{\ell}(T))=0$ if and
only if $\E(f|\mathcal{Z}_\ell(T_t))=0$ for $\sigma$-a.e. $t\in[0,1]$.

We note that for ergodic systems the factor $\cZ_0=\mathcal{I}$ is
trivial and  $\cZ_1=\mathcal{K}$. The factors $\cZ_\ell$   are of
particular interest because they can be used to study the
limiting behavior  in $L^2$ of some multiple ergodic averages.
\begin{theorem}[{\bf Leibman~\cite{L1}}]\label{hkl}
Let $p_1,p_2,\ldots,p_s$ be a family of essentially distinct polynomials $($$p_i$, and 
$p_i-p_j$ for $i\neq j$ are nonconstant$)$ with integer coefficients. There exists a nonnegative integer $\ell=\ell(p_1,p_2,\ldots,p_s)$ with the following property:
If $(X,\cB,\mu,T)$ is a system and $f_1,f_2,\ldots,f_s\in L^\infty(X)$, then  the limit
$$
\lim_{N\to+\infty}\frac1N\sum_{n=1}^N T^{p_1(n)}f_1\cdot T^{p_2(n)}f_2\cdot\ldots\cdot T^{p_s(n)}f_s
$$
exists in $L^2(\mu)$; and it is  equal to zero as long as one of the functions $f_i$ is orthogonal to the factor $\mathcal{Z}_\ell(T)$.
\end{theorem}
(We say that $\mathcal{Z}_\ell(T)$ is a {\it characteristic factor} associated with $p_1,p_2,\ldots,p_s$ when this last fact is true.)

Here are some examples that will be used later:
\begin{itemize}\item[(i)]
\cite{HK1} If $p_i(n)=in,\ 1\leq i\leq s$, then $\ell(p_1,p_2,\ldots,p_s)=s-1$.\item[(ii)]
\cite{F} More generally, if $p$ is a nonconstant integer polynomial and if $p_i(n)=ip(n),\ 1\leq i\leq s$, then $\ell(p_1,p_2,\ldots,p_s)=s-1$.
\item[(iii)]
\cite{FK2}  If the polynomials $p_1,p_2,\ldots,p_s$ are linearly independent and have zero constant term, then $\mathcal{K}_{rat}(T)$ is a characteristic factor.
\end{itemize}

\begin{proposition}\label{dc1}
Let $(X,\cB,\mu,T)$ be a system,  $p_1,p_2,\ldots,p_s$ be a family of  essentially distinct polynomials with integer coefficients, and let  $\mathcal{Z}_\ell(T)$ be a characteristic factor associated with this family.
 If  $f_0,f_1,\ldots,f_s\in L^\infty(X)$ and  one of the functions $f_i$ is orthogonal to the factor $\mathcal{Z}_{\ell+1}(T)$, then
$$
\text{D-}\!\!\lim_{n\to\infty} \int f_0\cdot
T^{p_1(n)}f_1\cdot\ldots\cdot T^{p_s(n)}f_s \ d\mu =0.
$$
\end{proposition}

\begin{proof}
As we mentioned before, if $f_i$ is orthogonal to the factor $\mathcal{Z}_{\ell+1}(T)$, then  $f_i\otimes \overline{f_i}$ is orthogonal to the factor $\mathcal{Z}_{\ell}(T\times T)$. By remark $(iii)$ following Theorem~\ref{hkl},  the averages
$$
\frac1N\sum_{n=1}^N\int f_0(x)
\cdot \overline{f_0}(y)\cdot f_1(T^{p_1(n)}x)\cdot \overline{f_1}(T^{p_1(n)}y)\cdot \ldots \cdot f_s(T^{p_s(n)}x)\cdot \overline{f_s}(T^{p_s(n)}y)\ d\mu(x)d\mu(y)
$$
converge to zero. This gives that
$$
\lim_{N\to+\infty}\frac1N\sum_{n=1}^N\left|\int f_0\cdot
T^{p_1(n)}f_1\cdot\ldots\cdot T^{p_s(n)}f_s \ d\mu\right|^2=0,
$$
and proves the announced convergence in density.
\end{proof}

Similarly, we have the following result:
\begin{proposition}\label{dc2}
Let $p_1,p_2,\ldots,p_s$ be a family of linearly independent integer polynomials  with zero constant term.
Let $(X,\cB,\mu,T)$ be a system and $f_0,f_1,\ldots,f_s\in L^\infty(X)$.
If one of the functions $f_i$ is orthogonal to the Kronecker factor $\mathcal{Z}_{1}(T)$, then
$$
\text{D-}\!\!\lim_{n\to\infty} \int f_0\cdot
T^{p_1(n)}f_1\cdot\ldots\cdot T^{p_s(n)}f_s \ d\mu =0.
$$
\end{proposition}
\begin{proof} By remark $(iii)$ following Theorem~\ref{hkl}, a characteristic factor associated with the family of polynomials $p_1,p_2,\ldots,p_s$  is $\mathcal{K}_{rat}$. Furthermore, it is well known that if  $f_i$ is orthogonal to the factor $\mathcal{Z}_{1}(T)$, then  $f_i\otimes \overline{f_i}$ is orthogonal to the factor $\mathcal{K}_{rat}(T\times T)$. The same argument as in the preceding proof applies.
\end{proof}

\subsection{Nilsystems}\label{SS:nilsystems}
We will now define a class of systems of purely algebraic structure
 that will be crucial for our study.
 Given a topological group $G$, we denote the identity
element by $e$ and
  we let $G_0$ denote the connected component of $e$.
  If $A, B\subset G$, then $[A,B]$ is
defined to be the subgroup generated by elements of the form
$\{[a,b]:a\in A, b\in B\}$ where $[a,b]=ab a^{-1}b^{-1}$. We
define the commutator subgroups recursively by $G_1=G$ and
 $G_{\ell+1}=[G, G_{\ell}]$. A group
$G$ is said to be {\it $\ell$-step nilpotent} if its $(\ell+1)$
commutator $G_{\ell+1}$ is trivial. If $G$ is an $\ell$-step nilpotent
Lie group and $\gG$ is a discrete cocompact subgroup, then the
compact space $X = G/\gG$ is said to be an {\it $\ell$-step
nilmanifold}.  The group $G$ acts on $G/\gG$ by left translation
where the translation by a fixed element $a\in G$ is given by
$T_{a}(g\gG) = (ag) \gG$. Let $m$ denote the unique probability
measure on $X$ that is invariant under the action of $G$ by left
translations (called the {\it Haar measure}) and let $\G/\gG$
denote the Borel $\sigma$-algebra of $G/\gG$. Fixing an element
$a\in G$, we call the system $(G/\gG, \G/\gG, m, T_{a})$ an {\it
$\ell$-step nilsystem}.

Nilsystems play a central role in our study  because they provide a sufficient class
 for verifying several multiple recurrence results for
general measure preserving systems. In fact when one deals with ``polynomial recurrence" this  is  usually a consequence of Theorem~\ref{hkl} and the following result of Host and Kra (a closely related result was subsequently proved by Ziegler (\cite{Z})):
 \begin{theorem}[{\bf Host \& Kra~\cite{HK1}}]\label{T:HoKra}
Let $(X,\mathcal{B},\mu,T)$ be a system. Then, for every $\ell\in\N$, a.e. ergodic component of the factor $\mathcal{Z}_\ell(T)$ is
an inverse limit of  $\ell$-step nilsystems.
\end{theorem}

Fundamental properties of nilsystems, related to our discussion,
were studied in \cite{AGH}, \cite{Pa}, \cite{Les}, and \cite{L}.
Below we summarize some facts that we shall  use, all the proofs
can be found in \cite{L}.

If $H$ is a closed subgroup of $G$ and $x\in X$, then $Hx$ may not be a  closed subset of $X$ (take $X=\R/\Z$, $x=\Z$, and
$H=\{k\sqrt{2}\colon k\in \Z\}$), but if $Hx$ is closed in $X$,  then the compact set $Hx$  can be given the structure of a nilmanifold. More precisely, if $x=g\Gamma$ and $Hx$ is closed,
 we have   $Hx\simeq
H/\Delta$ where $\Delta=H\cap g\Gamma g^{-1}$, and $h\mapsto hg\Gamma$ induces the isomorphism from $H/\Delta$ onto $Hx$. We call any such set a {\it
sub-nilmanifold} of $X$.

Let $(X=G/\Gamma,\mathcal{G}/\Gamma,m,T_a)$ be an ergodic
nilsystem. The subgroup  $\langle G_0,a\rangle$ projects to an open subset of
$X$ that is invariant under $a$. By ergodicity this projection
equals $X$. Hence, $X=\langle G_0,a\rangle/\Gamma'$ where $\Gamma'=\Gamma\cap
\langle G_0,a\rangle$. Using this representation of $X$  for  ergodic
nilsystems, we have that $G$ is generated by $G_0$  and  $a$.
From now on, when we work with an ergodic nilsystem, we will freely
assume that this extra  hypothesis  is satisfied, and so for example
we can assume that the commutator subgroups $G_\ell$ are connected for $\ell\geq 2$
(\cite{Les}).

\subsection{Uniform distribution properties in nilmanifolds}
If $G$ is a nilpotent Lie group,  $a_1,\ldots, a_s \in G$, and
$p_1, \ldots, p_s $ are integer polynomials
$\mathbb{N}\to\mathbb{Z}$, then a sequence of the form
 $g(n) =a_1^{p_1(n)}a_2^{p_2(n)}\cdots$ $a_s^{p_s(n)}$
  is called a {\it polynomial sequence} in $G$. In the sequel we need to establish
  various
uniform distribution properties of polynomial sequences on
nilmanifolds. The next result will simplify our task:
\begin{theorem}[{\bf Leibman~\cite{L}}]\label{T:L}
Let  $X = G/\gG$ be a  nilmanifold and $g(n)$ be a polynomial
sequence in $G$. Define $Z=G/([G_0,G_0]\Gamma)$ and let $\pi_Z\colon
X\to Z$ be the natural projection. For every $x\in X$:

(i) There exist $x_i\in X$ and connected sub-nilmanifolds $Y_i=Hx_i$ of $X$ (not necessarily distinct),
$1\leq i\leq t$, where
  $H$  is a closed subgroup of $G$ (depending on
$x$),   such that
$Y=\overline{\{g(n)x\colon n\in\mathbb{N}\}}=\bigcup_{i=1}^t Y_i$, and
for $i=1,\ldots,t$ the sequence $(g(tn+i))_{n\in\N}$ is
uniformly distributed in $Y_i$. If $Y$ is connected, then $t=1$.

 (ii) If $X$ is connected, then the sequence $(g(n)x)_{n\in\mathbb{N}}$ is dense in
$X$ if and only if it is  uniformly distributed in  $X$. Moreover,
$(g(n)x)_{n\in\mathbb{N}}$ is dense in $X$ if and only if
$(g(n)\pi_Z(x))_{n\in\mathbb{N}}$ is dense in $Z$.
\end{theorem}
We remark that the groups $G_0$ and $[G_0,G_0]$ are  normal
subgroups of $G$. The group   $G/[G_0,G_0]$ has the additional
property that the connected component of its identity element is
Abelian. This forces every rotation on the nilmanifold
$Z=G/([G_0,G_0]\Gamma)$ to have very special structure. More
precisely, a map $T\colon G\to G$ is said
to be {\it affine} if $T(g) = bA(g)$ for an endomorphism $A$ of
$G$ and some $b\in G$. Let $\ell\in\N$; the endomorphism $A$, or the affine
transformation $T$, is said to be {\it $\ell$-step unipotent} if $(A-{\text Id})^\ell=0$.
\begin{theorem}[{\bf F. \& Kra~\cite{FK}}]\label{T:FK}
Let $X=G/\gG$ be an $\ell$-step  connected nilmanifold with the Haar measure $m$
such that $G_0$ is Abelian and $a\in G$. 

Then the nilsystem $(X,
\mathcal{G}/\Gamma, m, T_a)$  is isomorphic to an $\ell$-step
unipotent affine transformation on some finite dimensional torus $\T^d$
with the Haar measure. Furthermore, the conjugation can be taken to be continuous.\end{theorem}


\section{Good and bad powers for sets of
$\ell$-recurrence}\label{S:multiple1} In this section we will prove
Theorem~\ref{T:B'}, but  before delving into the proof  let us motivate a
bit the choice of the set $R$. Suppose we just want to construct a
set $R$ that is bad for double recurrence but the set of squares of its elements
$R^2$ is good for double recurrence. In view of
Proposition~\ref{P:powers} below, it makes sense to take
$R=\big\{n\in\mathbb{N}\colon\{p(n)\beta\}\in [1/4,3/4]\big\}$ for some
quadratic polynomial $p$, this way we guarantee that $R$ will be
bad for double recurrence. It remains to choose the polynomial $p$ such that the set $R^2$ is good for double recurrence.  The obvious choice $p(n)=n^2$ will not
work since then $R^2$ will not even be good for single recurrence.
But the choice $p(n)=n^2+n$ will do the job (for any irrational
$\beta\in \R$) and this will be formally shown using
Proposition~\ref{P:ergodic2} below.

\subsection{Proof of the main theorem  modulo a multiple ergodic theorem}
We shall first establish Theorem~\ref{T:B'} modulo an ergodic theorem that
we will prove in subsequent sections. We need one preliminary
result that  was proved in \cite{FLW} in the special case
where the polynomial $p$ is a monomial. A very similar argument gives
the following more general result:
\begin{proposition}\label{P:powers} Let $R$ be a set of $\ell$-recurrence
and $p$ be an integer polynomial with zero constant term and
$\deg{p}\leq \ell$. For every $\alpha\in\mathbb{R}$ and
$\varepsilon>0$ there exists  $r\in R$ such that $\{p(r)\alpha\}\in
[0,\varepsilon]\cup [1-\varepsilon,1)$.
\end{proposition}
\begin{example*}
Let us illustrate how one proves Proposition~\ref{P:powers} in the case where
$\ell=2$ and $p(n)=n^2+n$.

Let $0<\varepsilon<1$, $\alpha\in\R$ and set  $\alpha'=\alpha/2$. It is well known that the set
$$
 \Lambda=\big\{n\in\N\colon \{ n\alpha \} \text{ and } \{n^2 \alpha'\}\in
 [0,\varepsilon/4]\big\}
$$
has positive density. Since $R$ is a set of double recurrence and
$\Lambda$ has positive density,   there exist $m\in \N$ and nonzero $r\in
R$ such that $m,m+r,m+2r \in \Lambda$. Then if
$$
A=\{m\alpha\}, B=\{(m+r)\alpha\}, C=\{m^2\alpha'\}, D=\{(m+r)^2\alpha'\}, E=\{(m+2r)^2\alpha'\}
$$
we have that $A,B,C,D,E\in [0,\varepsilon/4]$. Since
$$
B-A= r\alpha \pmod{1}, \quad C+E-2D=r^2\alpha \pmod{1},
$$
we have that  $$\{r\alpha\},\  \{r^2\alpha\} \in [0,\varepsilon/2]\cup [1-\varepsilon/2, 1).$$
It follows that $\{(r^2+r)\alpha\} \in [0,\varepsilon]\cup [1-\varepsilon, 1)$.
\end{example*}

We will also use a multiple ergodic theorem, its proof uses deeper
results from ergodic theory and will be given in the next section.
\begin{proposition}\label{P:ergodic2}
Let $(X,\mathcal{B},\mu, T)$ be a  system, $f_1,\ldots, f_\ell\in
L^\infty(\mu)$, $h_1,\ldots,h_s\colon \mathbb{T}\to\mathbb{C}$
be Riemann integrable functions,  and $\beta$ be an irrational
number. If $k,k_1,\ldots,k_s$ are distinct positive integers,
 then
\begin{align}\label{E:main}
\lim_{N\to\infty}\frac{1}{N}\sum_{n=1}^N &
h_1\big((n^{{\ell k_1}}+n^{k_1})\beta\big)\cdot\ldots\cdot
h_s\big((n^{\ell k_s}+n^{k_s})\beta\big)\cdot T^{n^{k}}f_1\cdot \ldots
\cdot
T^{\ell  n^k}f_\ell=\\
 \notag  &\int h_1 \ dt \cdot\ldots \cdot \int h_s \ dt\cdot
 \lim_{N\to\infty}\frac{1}{N}\sum_{n=1}^N T^{n^k}f_1\cdot \ldots \cdot
T^{ \ell n^k}f_\ell,
\end{align}
where the convergence takes place in $L^2(\mu)$.
\end{proposition}
\begin{remarks}
$(i)$ The second limit appearing in the statement exists by  \cite{L1} (or \cite{HK2}). To see that the first limit exists, apply \cite{L1} for the averages
 $$
 \frac{1}{N}\sum_{n=1}^N
S^{n^{{\ell k_1}}+n^{k_1}}\tilde{h}_1 \cdot\ldots\cdot
S^{n^{{\ell k_s}}+n^{k_s}}\tilde{h}_s \cdot
S^{n^{k}}\tilde{f}_1\cdot \ldots \cdot
S^{\ell n^k}\tilde{f}_\ell
 $$
 where $S=T\times R_\beta$ acts on $X\times \T$, $\tilde{h}_i(x,t)=e(m_it), m_i\in\Z$, $\tilde{f}_i(x,t)=f_i(x)$,
 and then use an approximation argument.

$(ii)$ As it will become clear from the proof,  the integer polynomials
$n^{{ \ell k_1}}+n^{k_1},\ldots, n^{{\ell k_s}}+n^{k_s}$ can be replaced by
  any family of polynomials $p_1,\ldots,p_s$ with zero constant term  having the following
property: For  every nonzero polynomial $p$ of the form
$p(n)=c_1n^k+ c_2n^{2k}+...+c_\ell n^{\ell k}$ the set $\{p_1,...,p_s,p\}$
is linearly independent.
\end{remarks}
We will also need an extension of Theorem~\ref{T:uniform}. In
order to not interrupt our discussion we state  and prove the
result needed  in the Appendix.
\begin{proof}[Proof of Theorem~\ref{T:B'}]
If $G=\N$, then by the polynomial extension of Szemer\'edi's theorem (\cite{BL}) $R=\N$ works. If $G\neq \N$, then  the complement $B$ of $G$ is nonempty. Let us first consider the case where $B$ is finite, say
$B=\{b_1,\ldots,b_s\}$, for some $b_i\in\N$. Fix an
irrational number $\beta$  and for $k\in\N$ let $R_k=\big\{n\in\N \colon
\{(n^{\ell k}+n^{k})\beta\} \in [1/4,3/4]\big\}$. We claim that the set
$R=R_{b_1}\cap \ldots \cap R_{b_s}$ has the advertised property.

Let  $b\in B$. If $p(n)=n^\ell +n$, then  the set  $p(R^b)$ is not
good for single recurrence for the rotation by $\beta$. It follows
from Proposition~\ref{P:powers} that the set $R^b$ is not good for
$\ell$-recurrence. On the other hand, let $g\in G$,
$(X,\mathcal{B},\mu,T)$ be a system and $A\in \mathcal{B}$ with
$\mu(A)>0$. We will show that $R^g$ is good for $\ell$-recurrence using
Proposition~\ref{P:ergodic2}. We set $h_i={\bf 1}_{[1/4,3/4]}$, $k_i=b_i$, for
$i=1,\ldots,s$, and $f_i={\bf 1}_A$ for $i=1,\ldots, \ell$ in
\eqref{E:main}, then multiply by $ {\bf 1}_A$ and integrate with
respect to $\mu$. We get
\begin{gather}
\notag\lim_{N\to\infty} \frac{1}{N}\sum_{1\leq n\leq N, n\in R}
\mu(A\cap T^{-n^g}A\cap\ldots\cap T^{-\ell n^g}A)=\\
\notag\frac{1}{2^s} \lim_{N\to\infty}\frac{1}{N}\sum_{n=1}^N
\mu(A\cap T^{-n^g}A\cap \ldots \cap T^{-\ell n^g}A).
\end{gather}
Using  the polynomial extension of Szemer\'edi's theorem (\cite{BL}) we get that the last limit is positive, proving that
$R^g$ is a set of $\ell$-recurrence.

The case where the set of bad powers $B$ is infinite is treated as in the proof of
Theorem~\ref{T:A'}. We  use the finite case and
Theorem~\ref{T:uniform-gen}.
\end{proof}
\subsection{Proof of the multiple ergodic theorem}\label{SS:proof}
In order to prove Proposition~\ref{P:ergodic2} we will use a
``reduction to affine" technique. The first step is a reduction to nil-systems and is based on Theorems~\ref{hkl} and \ref{T:HoKra}. Next we show that it suffices to verify the
result for a particular class of nilsystems, namely, for unipotent
affine transformations on finite dimensional tori (the
reduction is done in the course of proving Lemma~\ref{L:main}).
Lastly, we verify the result for such transformations (the main
ingredient is Lemma~\ref{L:affine}).

 We now execute our plan.
\begin{lemma}\label{L:affine}
Let $T$ be an  $\ell$-step unipotent affine transformation acting on
some finite dimensional torus $\mathbb{T}^d$, $\beta$ be an
irrational number, and $k,k_1,\ldots,k_s\in\N$ be distinct for
some $s\in \N$. For every  $x\in \T^d$ and functions $f\in
C(\T^d)$, $h_1,\ldots, h_s\in C(\T)$, we have
\begin{align}\label{E:affine}
\lim_{N\to\infty}\frac{1}{N}\sum_{n=1}^N &
h_1\big((n^{{\ell  k_1}}+n^{k_1})\beta\big)\cdot\ldots\cdot
h_s\big((n^{\ell k_s}+n^{k_s})\beta\big)\cdot f(T^{n^{k}}x)=\\
 \notag  &\int h_1 \ dt \cdot\ldots \cdot \int h_s \ dt\cdot
 \lim_{N\to\infty}\frac{1}{N}\sum_{n=1}^N f(T^{n^k}x).
\end{align}
\end{lemma}
\begin{proof}
  Without loss of generality we can assume that
$k_1<k_2<\ldots<k_s$. Arguing as in the proof of \ref{P:ergodic1},
it suffices to verify \eqref{E:affine} when $f(x)=\chi(x)$
is a  character of $\mathbb{T}^d$ and $h_1(t)=e(c_1t),\ldots,
h_s(t)=e(c_st)$, where $c_1,\ldots, c_s\in \Z$. Equivalently, we need to show
that if one of the $c_i$'s is nonzero, then for every $x\in
\mathbb{T}^d$  we have
\begin{equation}\label{E:reduced}
\lim_{N\to\infty}\frac{1}{N}\sum_{n=1}^N
e\big((c_1(n^{{\ell k_1}}+n^{k_1})+\ldots+
c_s(n^{\ell k_s}+n^{k_s}))\beta\big)\cdot  \chi(T^{n^{k}}x)=0.
\end{equation}
Let $i_{min}$ be the minimum $i\in \{1,\ldots, s\}$ such that
$c_i\neq 0$ and $i_{max}$ be the maximum $i\in \{1,\ldots, s\}$
such that $c_i\neq 0$.  Since the affine
transformation  $T$ is $\ell$-step unipotent we get that the coordinates of $T^nx$ are polynomials in
$n$ of degree at most  $\ell$. Indeed, if $Tx=Ax+b$, for some $b\in \T^d$ and endomorphism $A$ of $\T^d$ that satisfies $(A-I)^{\ell}=0$, then for $n\geq \ell$ we have
$$T^nx=\sum_{k=0}^{\ell-1} \binom{n}{k}(A-I)^kx+\sum_{k=0}^{\ell-1} \binom{n}{k+1}(A-I)^kb.$$
 It follows that
$$
\chi(T^{n^k}x)=e(q(n^k))
$$
for some real valued polynomial $q$ (depending on $x$) with $\deg q\leq \ell $. We denote $p(n)=q(n^k)$.
Consider two cases:

{\bf Case 1}. Suppose that $k_{i_{min}}<k$.  Notice that the non-constant terms of the polynomial
 $p$
have degree greater  or equal to $k$.
Therefore, in \eqref{E:reduced} we are averaging a sequence of the
form $e(P(n))$ for some polynomial $P$ that has a nonconstant
irrational coefficient (namely the coefficient of
$n^{k_{i_{min}}}$). By Weyl's uniform distribution theorem (\cite{We}) this
average must converge to zero.

{\bf Case 2}. Suppose that $k_{i_{min}}>k$. Since $\deg p\leq \ell k$,
  in \eqref{E:reduced} we are again averaging a
sequence of the form $e(P(n))$ for some polynomial $P$ that has
a nonconstant irrational coefficient (namely the coefficient of
$n^{\ell k_{i_{max}}}$). We conclude by Weyl's uniform distribution
theorem that this average must converge to zero.
\end{proof}

\begin{lemma}\label{L:main}
Let $(X=G/\Gamma, \mathcal{G}/\Gamma, m, T_a)$ be an $\ell $-step
nilsystem, $x\in X$, $\beta$ be an irrational number, and $k,k_1,\ldots,
k_s$ be distinct positive integers for some $s\in \N$.
Let
$g(n)=\big((n^{\ell k_1}+n^{k_1})\beta,\ldots,
(n^{\ell k_s}+n^{k_s})\beta, a^{n^k}x\big)$. The sequence
$(g(n))_{n\in \N} $ has an asymptotic distribution of the form $\lambda\otimes\rho$ where $\lambda$ is the Lebesgue measure on $\T^s$ and $\rho$ is the asymptotic distribution of the sequence $\big(a^{n^k}x\big)_{n\in \N}$ in $X$.\end{lemma}
\begin{proof}
Let us denote $Y=\overline{\{a^{n^k}x\colon n\in\N\}}$. Suppose first
that  the set  $Y$ is connected. It  follows from part $(i)$ of
Theorem~\ref{T:L} and the discussion in
Section~\ref{SS:nilsystems} that $Y$ is isomorphic to a connected
sub-nilmanifold $H/\Delta$ of $X$, so we can assume that $Y=H/\Delta$.


  We need to show that
$(g(n))_{n\in \N}$ is uniformly distributed on the nilmanifold
$\mathbb{T}^s\times Y$. By part $(ii)$ of Theorem~\ref{T:L} it
suffices to show that  the sequence
$\big(\big((n^{\ell k_1}+n^{k_1})\beta,\ldots, (n^{\ell k_s}+n^{k_s})\beta,
a^{n^k}\pi_Z(x)\big)\big)_{n\in\N}$ is uniformly distributed on
$\mathbb{T}^s\times Z$ where $Z=H/([H_0,H_0]\Delta)$ and $\pi_Z
\colon Y\to Z$ is the natural projection.  Substituting
$H/[H_0,H_0]$ for $H$ we can assume that $Z=H/\Delta$ where $H_0$ is Abelian.  Since
$Z$ is connected and $H_0$ is Abelian, by Theorem~\ref{T:FK} we
can assume that $T_a$, acting on $Z$, is an $\ell$-step unipotent affine
transformation on some finite dimensional torus. In this case the
result follows from Lemma~\ref{L:affine}.

In the general case we argue as follows: By part $(i)$ of
Theorem~\ref{T:L} we have $Y=\bigcup_{i=1}^t Y_i$ where $Y_i$ are
connected subnilmanifolds of $X$ such that
$Y_i=\overline{\{a^{(tn+i)^k}x\colon n\in\N\}}$ for $i=0,\ldots,t-1$.
Applying the previous argument (coupled with the analogous version
of Lemma~\ref{L:affine}) we get that
 for $i=0,\ldots,t-1$ the sequence $(g(tn+i))_{n\in\N}$
is uniformly distributed on the set  $ \mathbb{T}^s\times Y_i$.
This gives  the announced result with $\rho$ being the arithmetic mean of the uniform probabilities on the $Y_i$'s.
\end{proof}
\begin{proof}[{\bf Proof of Proposition~\ref{P:ergodic2}}]
First notice that using an ergodic decomposition argument we can
assume that the system is ergodic.
 Since by \cite{L1}  both limits in \eqref{E:main}
exist
 it suffices to show that
identity \eqref{E:main} holds weakly, that means, for $f_0,\ldots,
f_\ell \in L^\infty(\mu)$, and Riemann integrable functions
$h_1,\ldots,h_s\colon \mathbb{T}\to\mathbb{C}$ we have
\begin{align}\label{E:main'}
\lim_{N\to\infty}\frac{1}{N}\sum_{n=1}^N &
h_1\big((n^{{\ell k_1}}+n^{k_1})\beta\big)\cdot\ldots\cdot
h_s\big((n^{\ell k_s}+n^{k_s})\beta\big)\cdot\int f_0\cdot
T^{n^{k}}f_1\cdot \ldots \cdot
T^{\ell n^k}f_\ell \ d\mu=\\
 \notag  &\int h_1 \ dt \cdot\ldots \cdot \int h_s \ dt\cdot
 \lim_{N\to\infty}\frac{1}{N}\sum_{n=1}^N \int f_0\cdot T^{n^k}f_1\cdot \ldots \cdot
T^{\ell n^k}f_\ell\ d\mu,
\end{align}
If  $f_i\bot \mathcal{Z}_\ell $ for some $i\in \{0,1,\ldots,\ell \}$, then
by Proposition~\ref{dc1} and Example (ii) after Theorem~\ref{hkl}, we get that
$$
\text{D-}\!\lim_{n\to\infty} \int f_0\cdot T^{n^k}f_1\cdot \ldots \cdot
T^{\ell n^k}f_\ell\ d\mu,
$$
 and so both limits in \eqref{E:main'} are zero.
So we can assume that $f_i\in \mathcal{Z}_\ell $ for all $i\in
\{0,1,\ldots,\ell \}$. By \cite{HK1}, we know that the factor $\mathcal{Z}_\ell $ is
isomorphic to an inverse limit of $\ell$-step nilsystems. Moreover,
using a standard approximation argument we reduce our study to the case where the
system is an $\ell$-step nilsystem, say $(X=G/\Gamma,
\mathcal{G}/\Gamma, m, T_a)$ for some $a\in G$. Hence,
\eqref{E:main'} would follow if we show that for every $x\in X$
the sequence $\big(\big((n^{\ell k_1}+n^{k_1})\beta,\ldots,
(n^{\ell k_s}+n^{k_s})\beta, a^{n^k}x, \ldots,a^{\ell n^k}x\big)\big)_{n\in
\N}$ has an asymptotic distribution in $\T^s\times X^\ell $ of the form $\lambda\otimes\rho$, where $\lambda$ is the Lebesgue measure on $\T^s$ and $\rho$ is the asymptotic distribution of the sequence $\big(a^{n^k}x, \ldots,a^{\ell n^k}x\big)_{n\in \N}$ in $X^\ell $. But this follows from Lemma~\ref{L:main}
applied to the nilsystem induced by the rotation by $b=(a,
a^2,\ldots, a^\ell )$ on the $\ell$-step nilmanifold $X^\ell $ for the
diagonal point $(x,x,\ldots,x)\in X^\ell $.
\end{proof}

\subsection{Related results and questions}\label{SS:questions}
We discuss here some possible variations on Theorems~\ref{T:B} and \ref{T:B'}. We have
shown that if $G$ is a prescribed set of integers, then there
exists $R\subset \N$ such that the set  $R^g$ is good for
$\ell$-recurrence for all $g\in G$, and $R^b$ is bad for
$\ell$-recurrence for $b\in B=\N\setminus G$. A natural question is whether
it is possible to
strengthen this result and make $R^b$ have  bad $1$-recurrence
properties. In several cases this can be done, but there are some
limitations too. For example, if $R^g$ is good for $\ell$-recurrence,
then the set $R^{kg}$, $k=1,\ldots \ell $,  is good for $1$-recurrence
for all circle rotations (see Proposition~\ref{P:powers}).  We
can show that this is actually the only restriction.
\begin{theorem} Let $G\subset \N$ and $\ell\colon G\to \N$. If
$B=\N\setminus G$, then the condition
$$
B\cap kG=\emptyset, \text{ for }  1\leq k \leq \ell (g)
$$
is necessary and sufficient for the existence of a set
$R\subset\N$ such that
\begin{itemize}\item for all $g\in G$, the
set $R^g$ is good for  $\ell(g)$-recurrence,\item for all $b\in B$,
the set $R^b$ is bad for recurrence for some circle rotation.
\end{itemize}
\end{theorem}

The proof is analogous to the proof of Theorem~\ref{T:B'} so we are just going to sketch it. If $B$ is finite, say  $B=\{b_1,b_2,\ldots,b_s\}$, we fix an irrational number $\beta$ and  define $R$ to be  the  intersection of  the sets  $\big\{n\in
\N\colon \{n^{b_i}\beta\}\in [1/4,3/4]\big\}$ for $i=1,\ldots,s$. For $b\in B$ the set $R^b$ is obviously bad for the $1$-recurrence for the rotation by $\beta$. To show that $R^g$ is good for
$\ell(g)$-recurrence when $g\in G$, we study the limiting behavior of the following multiple ergodic averages:
$$
\frac1N\sum_{n=1}^Nh_1(n^{b_1}\beta)\cdot \ldots\cdot
h_s(n^{b_s}\beta)\cdot f(T^{n^g}x)\cdot f(T^{2n^g}x)\cdot \ldots \cdot f(T^{\ell(g)n^g}x).
$$
We can establish an ergodic theorem analogous to
Proposition~\ref{P:ergodic2} using  a minor
modification of the argument used in Section~\ref{SS:proof}, showing that  the set $R^g$ is good for $\ell(g)$-recurrence.  The case where $B$ is infinite can be treated using the finite case and Theorem~\ref{T:uniform-gen}, much like it was done in the proof of Theorem~\ref{T:A'}.

As we remarked before, if $R$ is a set of $2$-recurrence, then $R^2$ is a set of recurrence for circle rotations.
The same method shows that it is actually a set of recurrence for
rotations on any multidimensional torus. But is $R^2$ necessarily
a set of $1$-recurrence?
\begin{question}
If $R\subset \Z$ is a set of $2$-recurrence, is it true that
$R^2=\{r^2\colon r\in R\}$ is a set of  $1$-recurrence?
\end{question}
Another   closely related question is the following (a similar question was asked in \cite{BGL}):
\begin{question}\label{Q:question2}
If $R\subset \Z$ is a set of $\ell$-recurrence for every $\ell\in \N$, does the same hold for
the set $R^2=\{r^2\colon r\in R\}$?
\end{question}

Theorem~\ref{T:B'} is a model for a variety of multiple recurrence
results one may attempt to prove. For example, the following question is
related to a
plausible generalization of Theorem~\ref{T:B'}:

\begin{question}  Let  $p_1,\ldots,p_\ell $ be a family of integer
polynomials with zero constant term and $G$ be a set of positive integers. Does there
exist  a set $R\subset \N$ such that: for  $k\in\N$, the set  $R^k=\{r^k\colon r\in R\}$ is good for  recurrence
along the sequence $(p_1(n),\ldots,p_\ell (n))$ if and only if $k\in
G$?
\end{question}
 We are unable to give a positive answer because we lack detailed
information about the limiting behavior of multiple ergodic
averages along general polynomial schemes.

\section{Good and bad powers for sets of
$(n^{a_1},\ldots,n^{a_\ell })$-recurrence}\label{S:multiple2}
We begin with a multiple ergodic theorem that will be used in the proof of Theorem~\ref{T:C}.
We remind the reader that the set  $A_\ell $ consists of all
$\ell$-tuples $(a_1,a_2,\ldots,a_\ell )\in \mathbb{N}^\ell $ such that
$a_1<a_2<\ldots<a_\ell $.

\begin{proposition}\label{P:ergodic3}
Let $(X,\mathcal{B},\mu, T)$ be an ergodic   system,
$f_1,\ldots, f_\ell \in L^\infty(\mu)$, $h_1,\ldots,h_s\colon
\mathbb{T}\to\mathbb{C}$ be Riemann integrable functions.
Furthermore let $(a_{1,1},\ldots,a_{1,\ell })$, $\ldots$,
$(a_{s,1},\ldots, a_{s,\ell })$, $(b_1,\ldots,b_\ell )$$\in A_\ell $ be
distinct vectors, and suppose that the real numbers $1,\alpha_1,\ldots,
\alpha_s$ are  rationally  independent. If $p_i(n)=
n^{a_{i,1}}+\ldots+n^{a_{i,\ell }}$ for $i=1,\ldots, s$, then
\begin{align}\label{E:main11}
\lim_{N\to\infty}\frac{1}{N}\sum_{n=1}^N &
h_1\big(p_1(n)\alpha_1\big)\cdot\ldots\cdot
h_s\big(p_s(n)\alpha_s\big)\cdot T^{n^{b_1}}f_1\cdot \ldots \cdot
T^{n^{b_\ell }}f_\ell =\\
 \notag  &\int h_1 \ dt \cdot\ldots \cdot \int h_s \ dt\cdot
 \lim_{N\to\infty}\frac{1}{N}\sum_{n=1}^N T^{n^{b_1}}f_1\cdot
\ldots \cdot T^{n^{b_\ell }}f_\ell ,
\end{align}
where the convergence takes place in $L^2(\mu)$.
\end{proposition}
\begin{remark}
Both limits exist by   \cite{L1} (to deal with the first limit see Remark $(i)$ following Proposition~\ref{P:ergodic2}).
\end{remark}
\begin{proof}
First notice that using an ergodic decomposition argument we can
assume that the system is ergodic.
 Since the limits in \eqref{E:main11}
exist, it suffices to show that identity \eqref{E:main11} holds
weakly, that means, for $f_0,f_1,\ldots, f_\ell\in L^\infty(\mu)$,
and Riemann integrable functions $h_1,\ldots,h_s\colon
\mathbb{T}\to\mathbb{C}$ we have
\begin{align}\label{E:main22}
\lim_{N\to\infty}\frac{1}{N}\sum_{n=1}^N &
h_1\big(p_1(n)\alpha_1\big)\cdot\ldots\cdot
h_s\big(p_s(n)\alpha_s\big)\cdot \int f_0 \cdot
T^{n^{b_1}}f_1\cdot \ldots \cdot
T^{n^{b_\ell }}f_\ell \ d\mu=\\
 \notag  &\int h_1 \ dt \cdot\ldots \cdot \int h_s \ dt\cdot
 \lim_{N\to\infty}\frac{1}{N}\sum_{n=1}^N \int f_0\cdot T^{n^{b_1}}f_1\cdot
\ldots \cdot T^{n^{b_\ell }}f_\ell \ d\mu.
\end{align}
If  $f_i\bot \mathcal{K}$ for some $i\in \{0,1,\ldots,\ell \}$, then by
Proposition~\ref{dc2} both limits in \eqref{E:main22} are zero. So we
can assume that $f_i\in \mathcal{K}$ for all $i\in
\{0,1,\ldots,\ell \}$. Every ergodic Kronecker system is isomorphic to a
rotation on a monothetic compact  Abelian group with the Haar measure, and
any such group is the inverse limit of groups of the form $\T^k\times \Z_d$ for some nonnegative integers $k,d$.  Hence, using
a standard approximation argument we can furthermore assume that
our system is a rotation on
$\mathbb{T}^k\times \Z_d$ with the Haar measure $m$, and also  that  $f_i(x)=\chi_i(x)$,
$i=0,1,\ldots,\ell $, for some characters
$\chi_0,\chi_1,\ldots,\chi_\ell $ of $\mathbb{T}^k\times \Z_d$, and
$h_i(t)=e(l_it)$, $i=1,\ldots,s$, for some
$l_1,\ldots,l_s\in\mathbb{Z}$. If $l_i=0$ for $i=1,\ldots, s$, then
\eqref{E:main22} is obvious. If this is not the case, without loss
of generality we can assume that  $l_1\neq 0$. Then the right hand
side of \eqref{E:main22} is zero. Furthermore, notice that the
integral $\int \chi_0 \cdot T^{n^{b_1}}\chi_1\cdot \ldots \cdot
T^{n^{b_\ell }}\chi_\ell \ dm$ is either zero or has the form
$e(n^{b_1}\beta_1+\ldots +n^{b_\ell }\beta_\ell )$ for some
$\beta_1,\ldots,\beta_\ell \in \mathbb{R}$. Keeping in mind  these two
facts, we see that  in order to verify \eqref{E:main22}  it
suffices to show that
\begin{equation}\label{E:sdf}
\lim_{N\to\infty}\frac{1}{N}\sum_{n=1}^N
e(l_1p_1(n)\alpha_1+\ldots+l_sp_s(n)\alpha_s+n^{b_1}\beta_1+\ldots+
n^{b_\ell }\beta_\ell )=0.
\end{equation}
Since  $(a_{1,1},\ldots,a_{1,\ell })$ and $(b_1,\ldots,b_\ell )$ are
distinct vectors in $A_\ell$, there exists $j\in \{1,\ldots, \ell \}$
such that $a_{1,j}\neq b_i$ for all $i\in \{1,\ldots,\ell\}$. Then
the coefficient $\gamma$ of $n^{a_{1,j}}$ in \eqref{E:sdf} is
equal to $l_1\alpha_1$ plus an integer linear combination of the numbers
$\alpha_2,\ldots,\alpha_s$. Since $l_1\neq 0$ and the numbers $1,
\alpha_1,\alpha_2,\ldots,\alpha_s$ are rationally independent it
follows that $\gamma$ is nonzero. By  Weyl's uniform distribution
theorem (\cite{We}) we conclude that  \eqref{E:sdf} is satisfied, completing the proof.
\end{proof}
We will also use  the following simple lemma:
\begin{lemma}\label{P:powers2} If $R$ is a set of recurrence
 along the sequence $(n^{a_1},\ldots,n^{a_\ell })$, then for every
 $\alpha\in \mathbb{R}$ and $\varepsilon>0$ there
 exists $r\in R$ such that
 $\{(r^{a_1}+\ldots+r^{a_\ell })\alpha\}\in [0,\varepsilon]\cup[1-\varepsilon,1)$.
\end{lemma}
\begin{proof}
Consider the
system induced by the rotation by $\alpha$ on $\mathbb{T}$ with
the Haar measure, and let $A=[0,\varepsilon/\ell ]$. Our assumption easily implies
that there exists $r\in R$ such that  $\{r^{a_1}\alpha\},\ldots, \{r^{a_\ell }\alpha\}\in
[0,\varepsilon/\ell ]\cup[1-\varepsilon/\ell ,1)$. The
result follows.
\end{proof}

\begin{proof}[Proof of Theorem~B] If $G=A_\ell$  then by the  polynomial extension of Szemer\'edi's theorem (\cite{BL}) we have that  $R=\N$ works.
So we can assume that  $G\neq A_\ell$. Suppose first that the nonempty set $B=A_\ell\setminus G$ is finite and consists
of the vectors $(b_{1,1},\ldots,b_{1,\ell })$, $\ldots$,
$(b_{s,1},\ldots, b_{s,\ell })$. Let
$$
R=\big\{n\in\N \colon \{p_1(n)\alpha_1\},\ldots,\{p_s(n)\alpha_s\}
\in [1/4,3/4]\big\}
$$
where $p_i(n)= n^{b_{i,1}}+\ldots+n^{b_{i,\ell }}$, $i=1,\ldots, s $,
and $1,\alpha_1,\ldots, \alpha_s\in \mathbb{R}$ are rationally
independent. We claim that $R$ is a set of recurrence along the
sequence $(n^{a_1},\ldots,n^{a_\ell })$ if and only if
$(a_1,\ldots,a_\ell )\in
 G$.

Let  $(b_{i,1},\ldots,b_{i,\ell })\in B$. Since
$\{p_i(n)\alpha_i\}\in [1/4,3/4]$ for every $n\in R$, we get by
Lemma~\ref{P:powers2} that $R$ is not a set of recurrence along the sequence
$(n^{b_{i,1}},\ldots,n^{b_{i,\ell }})$. We will now use Proposition~\ref{P:ergodic3} to show that
if
$(g_1,\ldots,g_\ell )\in  G$, then $R$ is  a set of recurrence along the sequence
$(n^{g_1},\ldots,n^{g_\ell})$. Let
  $(X,\mathcal{B},\mu,T)$ be a system, and
$A\in \mathcal{B}$ with $\mu(A)>0$. Set $h_i={\bf
1}_{[1/4,3/4]}$ for $i=1,\ldots,s$, and $f_i={\bf 1}_A$ for
$i=1,\ldots,\ell $ in \eqref{E:main11}, multiply by $ {\bf 1}_A$ and
integrate with respect to $\mu$.  We find that
\begin{gather}
\notag\lim_{N\to\infty} \frac{1}{N}\sum_{1\leq n\leq N, n\in R}
\mu(A\cap T^{-n^{g_1}}A\cap\ldots\cap T^{-n^{g_\ell }}A)=\\
\notag\frac{1}{2^s} \lim_{N\to\infty}\frac{1}{N}\sum_{n=1}^N
\mu(A\cap T^{-n^{g_1}}A\cap \ldots \cap T^{-n^{g_\ell }}A).
\end{gather}
By the polynomial extension of Szemer\'edi's theorem (\cite{BL}) we have that the last limit is positive, proving that
$R$ is a set of recurrence along the sequence
$(n^{g_1},\ldots,n^{g_\ell })$.

The case where $B$ is infinite is treated as in the proof of
 Theorem~\ref{T:A'}
 using the finite case and  Theorem~\ref{T:uniform-gen}.
\end{proof}

\section{Powers of sequences and sufficient conditions for $\ell$-recurrence}\label{S:conjecture}
Useful sufficient conditions for a given set of integers $R$  to be good
for single recurrence
 were given in \cite{KM}. It was shown there  that
 a set  $R\subset \N$ is good for  recurrence if
 for every $d\in\N$ there exists a  sequence $(r_{d,n})_{n\in\N}$
  with values in the set  $R_d=\{r\in R\colon d!|r\}$ such that
  the sequence $(r_{d,n}\alpha)_{n\in\N}$
  is uniformly distributed $\pmod{1}$
 for every irrational number $\alpha$. Our objective in this section is to
 discuss analogous sufficient conditions for higher order
 recurrence.

\subsection{Sufficient conditions for single
recurrence}
 We first prove a single recurrence result, similar to the one given in \cite{KM},
 that will serve as a prototype  for the higher order statement we have in mind. The
 argument is very similar to the one used by
 Furstenberg~(\cite{Fu2}) to prove that the set of squares is good
 for single recurrence.
\begin{theorem}[{\bf Kamae \& Mend\`es-France} \cite{KM}]\label{T:single}
Let $(r_n)_{n\in\N}$ be   sequence of  
integers that satisfies

\noindent $(i)$ The sequence  $(r_n \alpha)_{n\in \N}$ is uniformly
distributed in $\mathbb{T}$ for every irrational $\alpha$.

\noindent $(ii)$ The set $\{n\in \N \colon d|r_n\}$ has positive upper
density for every $d\in \N$.

\noindent Then the set  $R=\{r_1,r_2,\ldots \}$ is good for single recurrence.
\end{theorem}
\begin{proof}
We first use assumption $(i)$ to show the following: If
$(X,\mathcal{B},\mu,T)$ is a system and $f\in L^\infty(\mu)$ is
such that $\E(f|\mathcal{K}_{rat})=0$, then
\begin{equation}\label{E:22}
\lim_{N\to\infty}\frac{1}{N}\sum_{n=1}^N \int \bar{f} \cdot
T^{r_n}f\ d\mu= 0.
\end{equation}
Let $\sigma_f$ be the spectral measure of the function $f$ with
respect to our system. Then the limit in \eqref{E:22} is equal to
\begin{equation}\label{E:33}
\lim_{N\to\infty}\frac{1}{N}\sum_{n=1}^N \int_{[0,1)} e(r_nt)\
d\sigma_f(t)= \int_{[0,1)}\Big(
\lim_{N\to\infty}\frac{1}{N}\sum_{n=1}^N e(r_nt)\Big) \
d\sigma_f(t).
\end{equation}
 Since by assumption
$$
\lim_{N\to\infty}\frac{1}{N}\sum_{n=1}^N  e(r_nt)=0
$$
for $t$ irrational and the measure $\sigma_f$  has no rational
point masses (since $\E(f|\mathcal{K}_{rat})=0$),  the limit in
\eqref{E:33} is zero and \eqref{E:22} follows.

Next we use \eqref{E:22} and assumption $(ii)$ to finish the
proof. Let $A\in\mathcal{B}$ with $\mu(A)>0$. Setting $f={\bf
1}_A-\E({\bf 1}_A|\mathcal{K}_{rat})$ in \eqref{E:22} we find that
$$
\lim_{N\to\infty}\Big(\frac{1}{N}\sum_{n=1}^N \mu(A\cap
T^{-r_n}A)-\frac{1}{N}\sum_{n=1}^N \int \mathbb{E}({\bf
1}_A|\mathcal{K}_{rat})\cdot T^{r_n}\E({\bf
1}_A|\mathcal{K}_{rat}) \ d\mu\Big)=0.
$$
Hence, it suffices to show that
\begin{equation}\label{E:44}
\limsup_{N\to\infty} \frac{1}{N}\sum_{n=1}^N \int \mathbb{E}({\bf
1}_A|\mathcal{K}_{rat})\cdot T^{r_n}\E({\bf
1}_A|\mathcal{K}_{rat}) \ d\mu >0.
\end{equation}
To see this, let $\varepsilon>0$ (to be determined later),  and
choose $d\in \N$ such that
$$
\norm{\mathbb{E}({\bf 1}_A|\mathcal{K}_{d}) -\mathbb{E}({\bf
1}_A|\mathcal{K}_{rat})}_{L^2(\mu)}\leq \varepsilon.
$$
Let $R\cap d\N=\{t_1,t_2,\ldots \}$. Since by assumption the set
$\{n\in \N \colon d|r_n\}$ has positive upper density, in order to show
\eqref{E:44} it suffices to show that
\begin{equation}\label{E:55}
\liminf_{N\to\infty} \frac{1}{N}\sum_{n=1}^N \int \mathbb{E}({\bf
1}_A|\mathcal{K}_{rat})\cdot T^{t_n}\E({\bf
1}_A|\mathcal{K}_{rat}) \ d\mu >0.
\end{equation}
 Using
the triangle inequality twice we see that the limit in
\eqref{E:55} is greater or equal than
$$
\liminf_{N\to\infty} \frac{1}{N}\sum_{n=1}^N \int \mathbb{E}({\bf
1}_A|\mathcal{K}_d)\cdot T^{t_n}\E({\bf 1}_A|\mathcal{K}_{d}) \
d\mu-2\varepsilon.
$$
For $g\in \mathcal{K}_d$ and $n\in \N$ we have $T^{t_n}g=g$, hence
the last expression is equal to
$$
\int  \mathbb{E}({\bf 1}_A|\mathcal{K}_{d})^2 \ d\mu
-2\varepsilon.
$$
Finally, since $\int  \mathbb{E}({\bf 1}_A|\mathcal{K}_{d})^2 \
d\mu\geq \mu(A)^2$, we see that it suffices to choose $\varepsilon
<  \mu(A)^2/2$  in order to obtain \eqref{E:55}. This completes
the proof.
\end{proof}

Using standard estimates on exponential sums one can use
Theorem~\ref{T:single} to deduce some well known results (see \cite{Sa}), for
example that the set of squares $\{n^2\colon n\in\N \}$ and the set of
shifted primes $\{p-1\colon p \text{ prime}\}$ are good for single
recurrence. See~\cite{Ble} for a stronger version of the preceding
result and several other applications.

\subsection{A counterexample for double recurrence}\label{counterexample}
As stated in Proposition~\ref{P:powers}, double recurrence for a
set $R$ forces nontrivial single recurrence properties for sets of
the form $p(R)$ where $p$ is any quadratic integer polynomial with
zero constant term. So in order to extend Theorem~\ref{T:single}
to double recurrence, one is lead to consider uniform distribution
properties of ``quadratic nature" as  possible substitutes for
 condition $(i)$ of Theorem~\ref{T:single}. We first show why the term
``quadratic nature" cannot be characterized using standard
quadratic polynomials only.
\begin{theorem}\label{T:main2}
There exists a  sequence of  positive
integers $(r_n)_{n\in\N}$   that satisfies:

\noindent $(i)$ If either  $\gamma$ or $\delta$ is irrational, then
the sequence $(r_n^2 \gamma+r_n\delta)_{n\in \N}$ is uniformly
distributed in $\mathbb{T}$.

\noindent $(ii)$ The set $\{n\in \N \colon d|r_n\}$ has
  positive density for every $d\in \N$.

\noindent $(iii)$ The set $R=\{r_1,r_2,\ldots \}$ is not good for double recurrence.
\end{theorem}
In order to prove Theorem~\ref{T:main2} we need two preparatory
lemmas. The first one is proven in \cite{H} using van der Corput's
inequality and some elementary manipulations of the resulting
exponential sums. It shows that some simple generalized quadratic
sequence is asymptotically orthogonal to  standard quadratic
sequences.
\begin{lemma}\label{L:p}
Suppose that $1,\alpha, \beta$ are rationally independent real
numbers. Then for every real valued polynomial  $p$ with
$\deg{p}\leq 2$ we have
$$
\lim_{N\to\infty} \frac{1}{N}\sum_{n=1}^N
e([n\alpha]n\beta+p(n))=0.
$$
\end{lemma}

The next result has similar context with
Proposition~\ref{P:powers}, it strengthens the single recurrence
properties that we can deduce when we know that a set is  good for double
recurrence.
\begin{lemma}\label{L:1}
Suppose that $1, \alpha, \beta$ are rationally independent real
numbers. If $R$ is a set of double recurrence, then for every
$\varepsilon>0$ there exists nonzero $r\in R$ such that $\{[r\alpha]r\beta\}\in
[0,\varepsilon]\cup [1-\varepsilon, 1)$.
\end{lemma}
\begin{remark}
The conclusion  actually holds for every $\alpha,\beta \in
\mathbb{R}$ but we will not use this.
\end{remark}
\begin{proof}
Let $0<\varepsilon<1$ and $\beta'=\beta/4$. From Lemma~\ref{L:p}
and Weyl's criterion for uniform distribution it follows that the
sequence $([n\alpha]n\beta', n\beta')$ is uniformly distributed in
$\mathbb{T}^2$. As a result, the set
$$
 \Lambda=\big\{n\in\N\colon \{[n\alpha]n\beta'\} \text{ and } \{n\beta'\}\in
 [0,\varepsilon/4]\big\}
$$
has positive density. Since $R$ is a set of double recurrence and
$\Lambda$ has positive density,   there exist $m\in \N$ and nonzero $r\in
R$ such that $m,m+r,m+2r \in \Lambda$. Then
\begin{equation}\label{E:three}
\{m\beta'\}, \{(m+r)\beta'\}, \{(m+2r)\beta'\}\in
[0,\varepsilon/4],
\end{equation}
 and if
$$
A= \{[m\alpha]m\beta'\},\quad B=
\{[(m+r)\alpha](m+r)\beta'\},\quad C=
\{[(m+2r)\alpha](m+2r)\beta'\}
$$
we have that $A,B,C\in [0,\varepsilon/4]$. Using the identity
$[a+b]=[a]+[b]+{\bf 1}_{\{a\}+\{b\}>1}(a,b)$ we get that
\begin{align}\label{E:ABC}
A+C-2B=&[m\alpha]m\beta'+([m\alpha]+2[r\alpha]+e_1+e_2)(m+2r)\beta'- \\
 \notag
&2([m\alpha]+[r\alpha]+e_3)(m+r)\beta' \pmod{1}\\
\notag =&4[r\alpha]r\beta' +(e_1+e_2)(m+2r)\beta'
-2e_3(m+r)\beta'\pmod{1}
\end{align}
for some $e_1, e_2,e_3\in \{0,1\}$. Since $\{A+C-2B\}\in
[0,\varepsilon/2]\cup [1-\varepsilon/2, 1)$ and
$\{(e_1+e_2)(m+2r)\beta' -2e_3(m+r)\beta'\}\in
[0,\varepsilon/2]\cup [1-\varepsilon/2, )$ (by \eqref{E:three}),
equation \eqref{E:ABC} gives that $\{4[r\alpha]r\beta'\}\in
[0,\varepsilon]\cup [1-\varepsilon, 1)$. Keeping in mind that
$4\beta'=\beta$, the proof is complete.
\end{proof}

\begin{proof}[{\bf Proof of Theorem~\ref{T:main2}}]
Let $\alpha$ and $\beta$ be two real numbers such that $1$, $\alpha$ and $\beta$ are rationally independent. We
claim that the set $R=\big\{n\in\N\colon \{[n\alpha]n\beta\} \in
[1/4,3/4]\big\}$ has the advertised properties. First notice that
by Lemma~\ref{L:p} the sequence $([n\alpha]n\beta)_{n\in\N}$ is uniformly
distributed in $\mathbb{T}$, so we deduce that $d(R)=1/2$.

We  verify $(i)$. Since the sequence
$(n^2\gamma+n\delta)_{n\in\N}$ is uniformly distributed in $\T$,
using Weyl's criterion and Lemma~\ref{L:p} we
can  easily derive  that the sequence
$\big(([n\alpha]n\beta,n^2\gamma+n\delta)\big)_{n\in\N}$ is uniformly distributed
in $\T^2$. It follows  that for  nonzero integers $k$ we have
$$
\lim_{N\to\infty}\frac{1}{N}\sum_{n=1}^N {\bf
1}_{[1/4,3/4]}(\{[n\alpha]n\beta\})\cdot e\big(k(n^2\gamma+n\delta)\big)=0.
$$
Since $d(R)=1/2$ this gives for  nonzero integers $k$  that
$$
\lim_{N\to\infty}\frac{1}{N}\sum_{n=1}^N
e\big(k(r_n^2\gamma+r_n\delta)\big)=0.
$$
Hence, the sequence
$(r_n^2 \gamma+r_n\delta)_{n\in \N}$ is uniformly distributed in
$\mathbb{T}$.

 Next we
verify $(ii)$. Using  Lemma~\ref{L:p} we can show as in the proof
of $(i)$ that
$$
\lim_{N\to\infty}\frac{1}{N}\sum_{n=1}^N {\bf
1}_{[1/4,3/4]}(\{[n\alpha]n\beta\})\cdot e(np/q)=0
$$
for every noninteger rational number $p/q$. It follows that the set
$R$ is uniformly distributed in arithmetic progressions, a
statement  stronger than $(ii)$.

Finally we verify $(iii)$. By construction, for every $r\in R$ we
have $\{[r\alpha]r\beta\}\in [1/4,3/4]$. Hence, Lemma~\ref{L:1} shows that
$R$ is not a set of double recurrence, completing the proof.
\end{proof}

\subsection{Conjecture for $\ell$-recurrence}  The  discussion
in the previous section indicates that in order to give necessary conditions for double
recurrence, uniform distribution properties in $\T$ of generalized quadratic
sequences
 of the form $([r_n\alpha]r_n\beta)_{n\in\N}$  should also be
taken into account.  In the same manner,  higher order recurrence
forces us to look for uniform distribution properties involving
more complicated generalized polynomials, the general form of
which we will not attempt to spell out. It has become apparent in
recent years (for example in \cite{GT1}, \cite{GT2},
 a similar problem arises in the study of asymptotics of
$\ell$-term arithmetic progressions in the prime numbers) that a
more efficient way to encode all these conditions is to work in a
non-Abelian setup, and look into uniform distribution properties of
linear sequences on nilmanifolds.
 Using this language, we can formulate
 what we think is a natural generalization
   of the result of Kamae \& Mend\`es-France
   (Theorem~\ref{T:single}), and thus state some potentially  sufficient conditions for $\ell$-recurrence
     in a very condensed form.

      Firstly, we need
to  extend the notion
   of an irrational
   rotation on $\T$ to general connected nilmanifolds:  Given a connected nilmanifold $X=G/\Gamma$,
an {\it irrational nilrotation in $X$} is an
   element $a\in G$ such that the sequence $(a^n\Gamma)_{n\in\N}$
   is uniformly distributed on $X$. We remark that if $a\in  G$ is an irrational nilrotation,  then $a^d$ is
   irrational for every $d\in\N$.

\begin{conjecture1}
Let $(r_n)_{n\in\N}$ be a sequence of  positive
integers that satisfies:

 \noindent $(i)$ For every connected $\ell$-step  nilmanifold $X$ and
  every irrational nilrotation $a$ in $X$ the sequence
 $(a^{r_n}\Gamma)_{n \in\N}$ is uniformly distributed in $X$.

\noindent $(ii)$ The set $\{n\in\N \colon d| r_n\}$ has positive
upper density for every $d\in \N$.

\noindent Then the set $R=\{r_1,r_2,\ldots \}$ is good for $\ell$-recurrence.
\end{conjecture1}
Notice that for $\ell=1$ the conjecture is true since it  easily reduces to
Theorem~\ref{T:single}. It can be shown that  conditions $(i)$ and $(ii)$
  are satisfied for nonconstant polynomials with zero constant term
(see \cite{F}),  and recent work of Green and Tao (see
\cite{GT1} and \cite{GT2}) indicates that they are probably
satisfied for the set of  shifted primes $\{p-1\colon p\text{
prime}\}$ (the shift is needed only for condition $(ii)$).

For $\ell=2 $ the definition of multiplication on a $2$-step nilpotent
Lie group is simple enough to find it beneficial to rewrite condition
$(i)$ in coordinates, much like it is done in the Appendix B of
\cite{GT1}. One is then led to consider  uniform distribution
properties involving {\it generalized quadratics}, that means,
functions $p\colon\mathbb{R}\to \mathbb{R}$ of the form $
p(t)=\sum_{i=1}^k [\alpha_it]\beta_it+\gamma t^2+\delta t+c$ where
$\alpha_i,\beta_i,\gamma,\delta,c\in\mathbb{R}$ for $i=1,\ldots,k$.
We say that a generalized quadratic is {\it irrational} if the
sequence $(p(n))_{n\in\N}$ is uniformly distributed in $\T$.

\begin{conjecture2}
Let $(r_n)_{n\in\N}$ be  a  sequence of  positive
integers that satisfies:

\noindent $(i)$ The sequence $(p(r_n))_{n\in\N}$ is  uniformly
distributed in $\mathbb{T}$ for every irrational generalized
quadratic $p$.

\noindent $(ii)$ The set $\{n\in\N \colon d| r_n\}$ has
positive upper density for every $d\in \N$.

\noindent Then the set $R=\{r_1,r_2,\ldots \}$ is good for double recurrence.
\end{conjecture2}



\subsection{Proof of conjecture I for positive density sets}

To add credibility to the previous conjectures we will verify
Conjecture I in a special case.
\begin{theorem}\label{T:conjecture}
Conjecture I holds if $d(R)>0$.
\end{theorem}
\begin{remark}
Our argument can  actually be used to verify Conjecture I for  any sequence $(r_n)_{n\in\N}$ for which
the nilfactor $\mathcal{Z}_\ell$ turns out to  be
characteristic for the multiple ergodic averages related to the family  $\{r_n, 2r_n,\ldots, \ell r_n\}$.
\end{remark}

We first prove a lemma that will help us deal with systems that
 have nontrivial ``periodic part''.

\begin{lemma}\label{L:r}
Let $d\in\N$ and $(N_k)_{k\in\N}$ be an increasing sequence of positive integers. Suppose that  $(r_n)_{n\in\N}$ is a  sequence of  positive
integers that satisfies:

 \noindent $(i)$ For every connected $\ell$-step  nilmanifold $X$ and
  every irrational nilrotation $a$ in $X$ the sequence
 $(a^{r_n}\Gamma)_{n \in \N}$ is uniformly distributed in $X$ with respect to the intervals
 $[1,N_k]$.

\noindent $(ii)$ The set $I_d= \{n\in\N \colon d| r_n\}$ has positive
density with respect to the sequence of intervals $[1,N_k]$.\\
 Then the subsequence
 $(r_n)_{n\in I_d}$ also satisfies property $(i)$.
\end{lemma}
\begin{proof}
Let $m$ denote the Haar measure on $X$. It suffices to show
that for every connected $\ell$-step nilmanifold and irrational
nilrotation $a\in X$ we have
\begin{equation}\label{E:main0}
\lim_{N_k\to\infty}\frac{1}{|I_d\cap N_k|}\sum_{n\in I_d\cap [1,N_k]}
f(a^{r_n}\Gamma)=\int f\ dm
\end{equation}
for all $f\in C(X)$.

 We will prove \eqref{E:main0}
 by induction on $\ell$. For $\ell=1$ we are reduced to the Abelian case
 and the proof is easy (if not, look at the proof of the inductive step).
  Suppose that the
statement holds for all connected $(\ell-1)$-step nilmanifolds, and
let $X=G/\Gamma$ be a connected $\ell$-step nilmanifold and $a\in G$
be an irrational nilrotation in $X$.

 We
start with a reduction. Since $G$ is $\ell$-step nilpotent, the
subgroup $\Gamma_\ell=G_\ell
\cap\Gamma$ is normal in $G$. So
$G/\Gamma_\ell$ is a group and $X=(G/\Gamma_\ell)/(\Gamma/\Gamma_\ell)$.
Using this representation for $X$ we can assume that
$\Gamma_\ell=\{e\}$ and so $G_\ell$ is a compact Abelian Lie group.
Since $G_\ell$ is connected if $\ell\geq 2$ (see discussion at the end of  Section~\ref{SS:nilsystems}) we can further assume that
it is some finite dimensional torus.

By $ \mathcal{F}$ we denote the set of $g\in C(X)$ with the
following property: there exists a character
$\chi$ of the torus $G_\ell$ such that for every $b\in G_\ell$ we have  $g(bx)=\chi(b)\cdot
g(x)$ for every $x\in X$.  It follows from \cite{Les'} (see the
proof of Proposition on page 121) that linear combinations of
functions in $\mathcal{F}$ form a dense subset of $C(X)$. So it
suffices to verify \eqref{E:main0} for functions in $\mathcal{F}$.

Let $g\in \mathcal{F}$. If the character $\chi$ defined before is
trivial, then $g$ is $G_\ell$-invariant and so it factors through the
connected $(\ell-1)$-step nilmanifold $(G/G_\ell)/(\Gamma/\Gamma_\ell)$.
Applying the induction hypothesis we get that \eqref{E:main0}
holds for $g$ in place of $f$.

Suppose now that there exists a nontrivial character $\chi$ of
$G_\ell$ such that $g(bx)=\chi(b)\cdot g(x)$ for $x\in X$. Integrating this equation with respect to $x$ and using that
 $x\mapsto bx$ is measure preserving gives $\int g\ dm=0$. So
it suffices to show that the limit in \eqref{E:main0} is zero when $g$ takes the place of $f$.
Since $G_\ell$ is a torus there exists $e_d\in G_\ell$ such that $e_d^d={\bf 1}$ and
$e_d^j\neq {\bf 1}$ for $j=1,\ldots d-1$, where ${\bf 1}$ is the identity
element in $G_\ell$.

Since $a$ is an irrational nilrotation, for every $j\in \N$,
the sequence $(a^{j+nd}\Gamma)_{n\in \N}$ is uniformly distributed in $X$.
It follows that  for $j=1,\ldots,d$, the
nilrotations $b_j=a\cdot e_d^j$ are also irrational on $X$. So condition $(i)$ gives that
$$
\lim_{N_k\to\infty}\frac{1}{N_k}\sum_{n=1}^{N_k}
g(b_j^{r_n}x)=0
$$
for $j=1,\ldots,d$.
Averaging over $j$ we get
\begin{multline}\label{E:213}
0= \lim_{N_k\to\infty}\frac{1}{N_k}\sum_{n=1}^{N_k}
\Big(\frac{1}{d}\sum_{1\leq j\leq d}
g(b_j^{r_n}x)\Big)=\lim_{N_k\to\infty}\frac{1}{N_k}\sum_{n=1}^{N_k}
\Big(g(a^{r_n}\Gamma)\cdot \frac{1}{d}\sum_{1\leq j\leq d}
\chi(e_d^{jr_n})\Big)=
\\
\lim_{N_k\to\infty}\frac{1}{N_k}\sum_{n=1}^{N_k}
\Big(g(a^{r_n}\Gamma)\cdot {\bf 1}_{I_d}(n)\Big) =\lim_{N_k\to\infty}\frac{|I_d\cap [1,N_k]|}{N_k}
\cdot
\lim_{N_k\to\infty}\frac{1}{|I_d\cap [1,N_k]|}\sum_{n\in I_d\cap [1,N_k]}
g(a^{r_n}\Gamma).
\end{multline}
Note that we used the nontriviality of $\chi$  to justify the
third equality. Since  by condition $(ii)$ the density of the set $I_d$ with respect to the intervals $[1,N_k]$ is positive, \eqref{E:213} gives that
$$
\lim_{N_k\to\infty}\frac{1}{|I_d\cap [1,N_k]|}\sum_{n\in I_d\cap [1,N_k]}
g(a^{r_n}\Gamma)=0.
$$
So \eqref{E:main0}
is satisfied with $g$ in place of $f$. This completes the proof.
\end{proof}

\begin{proof}[Proof of Theorem~\ref{T:conjecture}]
Let $\ell\in \N$.
It suffices to show
that for every  system $(X,\mathcal{B},\mu,T)$
and $f\in L^\infty(\mu)$,  nonnegative and not a.e. zero, we have
\begin{equation}\label{positive}
\limsup_{N\to\infty} \frac{1}{N} \sum_{n=1}^{N} \int f\cdot
T^{r_n}f\cdot \ldots \cdot T^{\ell r_n}f\ d\mu>0.
\end{equation}

We start with some reductions.  Using  an ergodic decomposition
argument we can assume that the system is ergodic. Furthermore,
since $d(R)>0$, by example $(i)$ right after Theorem~\ref{hkl}  and Proposition~\ref{dc1}, the factor
$\mathcal{Z}_\ell$ is characteristic for the multiple ergodic
averages appearing in \eqref{positive}. Notice also that the
projection of $f$ onto $\mathcal{Z}_\ell$ is also nonnegative and not
a.e. zero. Since by Theorem~\ref{T:HoKra} the factor $\mathcal{Z}_\ell$ is an
inverse limit of $\ell$-step nilsystems, we have reduced the problem
to establishing \eqref{positive} for such systems. Moreover, an
argument completely analogous to that of Lemma $3.2$ in \cite{FuK}
shows that the positiveness property \eqref{positive} is preserved
by inverse limits.
Hence, we can further assume that the system is an ergodic
$\ell$-step nilsystem.

 In this case, by Proposition $7.2$ of \cite{BHK}, there exists an
ergodic $\ell$-step nilsystem $(Y=G/\Gamma, m, T_a)$ and a continuous
function $F$ on $Y$  such that
\begin{equation}\label{E:nilsequence}
\int f\cdot T^{n}f\cdot \ldots \cdot T^{\ell n}f\ d\mu= F(a^n\Gamma)
\end{equation}
for every $n\in\N$.

We can assume that for some $d\in \N$ we have  $Y=\Z_d\times Y_0$
where $Y_0$ is a connected nilmanifold
(the torsion part is cyclic since $Y$ admits an ergodic nilrotation).
By changing the first coordinate of $a$ to  $0\in \Z_d$, we get an
irrational  element $b$  of the connected nilmanifold $\{0\}\times Y_0$ that satisfies
$b^d=a^d$.

Let $I_d=\{n\in\N\colon d|r_n\}$. Since $\bar{d}(I_d)>0$, there exists  an increasing
sequence  of integers $N_k$  such that
$\lim_{N_k\to\infty} |I_d\cap[1,N_k]|/N_k>0$.
  So in order to establish \eqref{positive}, it
suffices to show that
\begin{equation}\label{positive1}
\limsup_{N_k\to\infty} \frac{1}{|I_d\cap[1,N_k]|} \sum_{n\in I_d\cap
[1,N_k]} \int f\cdot T^{r_n}f\cdot \ldots \cdot T^{\ell r_n}f\ d\mu>0.
\end{equation}
 Using  \eqref{E:nilsequence}, and noticing that $a^{r_n}=b^{r_n}$ for $n\in
 I_d$,
 we see that the last limit is equal to
$$
\limsup_{N_k\to\infty} \frac{1}{|I_d\cap[1,N_k]|} \sum_{n\in I_d\cap
[1,N_k]} F(a^{r_n}\Gamma)= \limsup_{N_k\to\infty}
\frac{1}{|I_d\cap[1,N_k]|} \sum_{n\in I_d\cap [1,N_k]}
F(b^{r_n}\Gamma).
$$
Since $b$ is an irrational nilrotation of   $\{0\}\times Y_0$, by
Lemma~\ref{L:r} the sequence $(b^{r_n}\Gamma)_{n \in I_d}$   is
uniformly distributed in $\{0\}\times Y_0$ with respect to the sequence of intervals $[1,N_k]$. Also,
since $b^d$ is  irrational, the sequence
$(b^{dn}\Gamma)_{n \in\N}$ is  uniformly distributed in
$\{0\}\times Y_0$.
 It follows that the last  $\limsup$  is  equal to
$$
\lim_{N_k\to\infty} \frac{1}{|I_d\cap[1,N_k]|} \sum_{n\in I_d\cap
[1,N_k]} F(b^{r_n}\Gamma)= \lim_{N\to\infty} \frac{1}{N}
\sum_{n=1}^{N} F(b^{dn}\Gamma).
$$
Using once again that $b^d=a^d$ and \eqref{E:nilsequence}, we find
that the last limit is equal to
$$ \lim_{N\to\infty} \frac{1}{N} \sum_{n=1}^{N}
F(a^{dn}\Gamma)=\lim_{N\to\infty} \frac{1}{N} \sum_{n=1}^{N}\int
f\cdot T^{dn}f\cdot \ldots \cdot T^{d\ell n}f\ d\mu
$$
which  is positive by Furstenberg's multiple recurrence
theorem \cite{Fu1}. This  completes the proof.
\end{proof}

 \section{Appendix: Uniformity for sets of multiple recurrence}
We establish some uniform estimates that were used in the proofs
of Theorems~B'~and~C. Related estimates were obtained in
\cite{BHRF}, but our result is more general and our proof shorter.
\begin{definition}
A sequence $(u_1(n),\ldots,u_\ell(n))_{n\in\N}$ with values in $\Z^\ell$ is a {\it good sequence for multiple recurrence} if, given any system $(X,\mathcal{B},\mu,T)$ and any $A\in\mathcal{B}$ with $\mu(A)>0$, there exists  $n\in\N$ such that
$$
\mu(A\cap T^{-u_1(n)}A\cap \ldots\cap T^{-u_\ell(n)}A)>0.
$$
\end{definition}
Using Furstenberg's correspondence principle it is not hard  to verify  that the sequence $(u_1(n),\ldots,u_\ell(n))_{n\in\N}$ is good for multiple recurrence if and only if for  any set $\Lambda$ of positive upper density in $\Z$, there exists $n\in\N$ such that
$$
|\Lambda\cap (\Lambda-u_1(n))\cap \ldots\cap(\Lambda-u_\ell(n))|>0.
$$

\begin{theorem}\label{T:uniform-gen}
Let $(u_1(n),\ldots,u_\ell(n))_{n\in\N}$ be a good sequence for multiple recurrence. Then

$(i)$ For every $\varepsilon>0$ there exist
$\delta=\delta(\varepsilon)>0$ and $N_0=N_0(\varepsilon)$, such that
for every $N\geq N_0$ and   integer set
$\Lambda\subset[-N,N]$ with $|\Lambda|\geq\varepsilon N$, we have
$$
|\Lambda\cap (\Lambda-u_1(n))\cap \ldots \cap(\Lambda-u_\ell(n))|\geq \delta
N
$$
for some $n\in [1,N_0]$.

$(ii)$ For every $\varepsilon>0$ there exist
$\gamma=\gamma(\varepsilon)>0$ and $N_1=N_1(\varepsilon)$, such
that for every  system $(X,\mathcal{B},\mu,T)$ and $A\in
\mathcal{B}$ with $\mu(A) \geq\varepsilon$, we have that
$$
\mu(A\cap T^{-u_1(n)}A\cap \ldots\cap T^{-u_\ell(n)}A)\geq\gamma
$$
for some $n\in [1,N_1]$.
\end{theorem}
\begin{proof}
$(i)$  The argument is similar to one used in \cite{FK2}. Suppose
that the result fails. Then there exist $\varepsilon_0>0$,
sequence $(\delta_m)_{m\in\N}$ of positive real numbers with
$\lim_{m\to\infty} \delta_m=0$, increasing integer sequences
$(K_m)_{m\in\N},(N_m)_{m\in\N},$ with $N_m\geq K_m$, and integer sets
$\Lambda_m\subset[-N_m,N_m]$, such that
\begin{equation}
\label{finite1} |\Lambda_m|\geq \varepsilon_0N_m
\end{equation}
and
\begin{equation}
\label{finite2}|
\Lambda_m\cap(\Lambda_m-u_1(n))\ldots \cap(\Lambda_m-u_\ell(n))|\leq
\delta_m N_m\
\end{equation}
for every $m\in\N$ and $n\in [1,K_m]$. In order to get a
contradiction, we will construct a measure preserving system with
bad recurrence properties.

 For $m\in\N$ set $\Lambda^0_m=\Z\setminus \Lambda_m$ and $\Lambda^1_m=\Lambda_m$.
  Using a
diagonal argument we can find a subsequence of $(N_m)_{m\in\N}$,
which for convenience we denote  again by $(N_m)_{m\in\N}$, such
that the limit
$$
\lim_{m\to\infty} \frac{|(\Lambda^{i_1}_m-n_1)\cap
(\Lambda^{i_2}_{m}-n_2)\cap \ldots\cap (\Lambda^{i_s}_m-n_s)\cap[-N_m,N_m] |}{2N_m}
$$
 exists for every $s\in\N$, $n_1,\ldots, n_s\in\bbZ$, and
$i_1,\ldots,i_s\in\{0,1\}$.

On the sequence space $(X=\{0,1\}^\bbZ,\mathcal{B})$, where
$\mathcal{B}$ is the  Borel $\sigma$-algebra,  we define a measure
$\mu$ on cylinder sets as follows:
\begin{multline*}
\mu(\{x_{n_1}=i_{1},x_{n_2}=i_{2},\ldots,x_{n_s}=i_{s}\})=\\
\lim_{m\to\infty} \frac{|(\Lambda^{i_1}_m-n_1)\cap
(\Lambda^{i_2}_{m}-n_2)\cap \ldots\cap (\Lambda^{i_s}_m-n_s)\cap
[-N_m,N_m] |}{2N_m}
\end{multline*}
where $n_1, n_2, \ldots, n_s\in\bbZ$, and $i_{1}, i_2, \ldots,
i_s\in\{0,1\}$. Since $\mu$ is finitely additive on the algebra $\mathcal{F}$ of finite unions of cylinder sets, it is easy to 
check that $\mu$ defines a premeasure on $\mathcal{F}$ (see \cite{Mc}, Theorem 3.2.4), and hence, by
Carath\'eodory's extension theorem (\cite{Ca}),  it extends to a probability
measure on $\mathcal{B}$.
Then the shift
transformation $T$ defined by
$$T\big((x_j)_{j\in\bbZ}\big)=(x_{j+1})_{j\in\bbZ}
$$
preserves the measure $\mu$ (since this holds for $\mu$ restricted to $\mathcal{F}$) and gives rise to the system
$(X,\mathcal{B},\mu,T)$. If $$ A=\{x\in X \colon x(0)=1\},$$ using the
definition of $\mu$ we see that
\begin{align*}
\mu(A\cap T^{-u_1(n)}A\cap \ldots \cap T^{-u_\ell(n)}A)&=
\mu(\{x_{0}=1,x_{u_1(n)}=1,\ldots,x_{u_\ell(n)}=1\})\\
\notag &=\lim_{m\to\infty} \frac{|\Lambda_m\cap
(\Lambda_{m}-u_1(n))\cap \ldots\cap (\Lambda_m-u_\ell(n))|}{2N_m},
\end{align*}
for every $n\in\N$. Combining this with \eqref{finite1} and
\eqref{finite2}, and remembering  that $\lim_{m\to\infty}
\delta_m=0$, we find that $\mu(A)\geq \varepsilon_0/2>0$ and
$$\mu(A\cap T^{-u_1(n)}A\cap \ldots
\cap T^{-u_\ell(n)}A)=0
$$
for all $n\in \N$. This contradicts the fact that the sequence $(u_1(n),\ldots,u_\ell(n))_{n\in\N}$ is good for multiple recurrence
and
completes the proof of $(i)$.

$(ii)$ The argument is similar to  one used in \cite{BHRF}. Let
$\varepsilon>0$ and $A\in\mathcal{B}$ with
$\mu(A)\geq\varepsilon$. Set $N_1=N_0(\varepsilon/2)$ where $N_0$
was defined in part $(i)$. For $x\in X$ let
$$
f(x)=\frac{1}{N_1}\sum_{m=1}^{N_1} {\bf 1}_A(T^mx).
$$
Since $\int f \ d\mu=\mu(A)\geq\varepsilon$ and $0\leq f(x)\leq 1$, if $B=\{x\in X\colon
f(x)\geq\varepsilon /2\}$ we have that
$$
\mu(B)+(1-\mu(B))\cdot\varepsilon/2\geq \varepsilon,
$$
which implies that $\mu(B)\geq \varepsilon/2$. Notice that for
$x\in B$ we have that
$$
|\{m\in [1,N_1]\colon T^mx\in A\}|\geq \frac{\varepsilon}{2} N_1.
$$
Letting $ \Lambda_x=\{m\in [1,N_1]\colon T^mx \in A\}$ we get by
part $(i)$ that
$$
 \Lambda_x\cap(\Lambda_x-u_1(n))\cap \ldots \cap
(\Lambda_x-u_\ell(n))\neq {\varnothing},
$$
which implies  that
\begin{equation}\label{Lx}
{\bf 1}_A(T^mx)\cdot {\bf 1}_A(T^{m+u_1(n)}x)\cdot\ldots \cdot
{\bf 1}_A(T^{m+u_{\ell}(n)}x)=1
\end{equation}
 for some $m\in [1,N_1]$ and $n\in [1, N_1]$. Since we  have
   $N_1$ choices for $m$ and $n$,  we can choose
$m_0,n_0\in [1,N_1]$, and  a set $C\subset B$ such that
$\mu(C)\geq \mu(B)/N_1^2$ and \eqref{Lx} holds for all $x\in C$.
We have that
$$
\mu(A\cap T^{-u_1(n_0)}A\cap \ldots \cap T^{-u_\ell(n_0)}A)  = \int
{\bf 1}_A(T^{m_0}x)\cdot {\bf 1}_A(T^{m_0+u_1(n_0)}x)\cdot\ldots
\cdot {\bf 1}_A(T^{m_0+u_\ell(n_0)}x) \ d\mu
$$
which  is greater than
$$
 \int {\bf 1}_B(x)  \cdot {\bf
1}_A(T^{m_0}x)\cdot {\bf 1}_A(T^{m_0+u_1(n_0)}x)\cdot\ldots \cdot
{\bf 1}_A(T^{m_0+u_\ell(n_0)}x) \ d\mu\geq \mu(C),
$$
where the last inequality is valid because  \eqref{Lx} holds for
all $x\in C$. We have thus verified $(ii)$ with
$\gamma=\mu(B)/N_1^2\geq \varepsilon/(2N_1^2) $.
\end{proof}

\end{document}